\documentclass{article}

\usepackage[margin=1in]{geometry} 
\usepackage{amsmath,amsthm,amssymb}
\usepackage{graphicx}
\usepackage{color}
\usepackage[dvipsnames]{xcolor}
\usepackage{subcaption}
\usepackage[export]{adjustbox}
\usepackage{multirow}
\usepackage{url}

\DeclareMathOperator*{\argmin}{\text{argmin}}

\newcommand{\rev}[1]{#1}

\newtheorem{remark}{Remark}

\begin{document}

\title{Deep Global Model Reduction Learning in Porous Media Flow Simulation} 

\author{
Siu Wun Cheung\thanks{Department of Mathematics, Texas A\&M University, College Station, TX 77843, USA (\texttt{tonycsw2905@math.tamu.edu})}
\and
Eric T. Chung\thanks{Department of Mathematics, The Chinese University of Hong Kong, Shatin, New Territories, Hong Kong SAR, China (\texttt{tschung@math.cuhk.edu.hk}) }
\and
Yalchin Efendiev\thanks{Department of Mathematics \& Institute for Scientific Computation (ISC), Texas A\&M University,
College Station, Texas, USA (\texttt{efendiev@math.tamu.edu})}
\and
Eduardo Gildin\thanks{Department of Petroleum Engineering, Texas A\&M University, College Station, TX 77843, USA (\texttt{egildin@tamu.edu})}
\and
Yating Wang\thanks{Department of Mathematics, Texas A\&M University, College Station, TX 77843, USA (\texttt{wytgloria@math.tamu.edu})}
\and
Jingyan Zhang\thanks{Department of Mathematics, Texas A\&M University, College Station, TX 77843, USA (\texttt{jingyanzhang@math.tamu.edu})}
}

\date{}

\maketitle

\begin{abstract}

In this paper, we combine deep learning concepts and some proper orthogonal decomposition (POD) model
reduction methods for predicting flow in heterogeneous porous
media. Nonlinear flow dynamics is studied, where the dynamics
is regarded as a multi-layer network. The solution at the current time
step is regarded as a multi-layer network of the solution at the initial
time and input parameters. As for input, we consider various sources, which
include source terms (well rates), permeability fields, and initial
conditions. We consider the flow dynamics, where the solution is known
at some locations and the data is integrated to the flow dynamics by modifying the 
reduced-order model. This approach allows modifying the reduced-order
formulation of the problem. Because of the small problem size,
limited observed data can be handled. We consider enriching the observed data
using the computational data in deep learning networks. The basis functions
of the global reduced order model
are selected such that the degrees of freedom represent the solution
at observation points. This way, we can avoid learning basis functions,
which can also be done using neural networks. We present numerical results,
where we consider channelized permeability fields, where the 
network is constructed for various channel configurations.
Our numerical results show that one can achieve a good approximation
using forward feed maps based on multi-layer networks.


\end{abstract}

\section{Introduction}
\label{sec:intro}

Mathematical models are widely used to describe the underlying physical process
in science and engineering disciplines. 
In many real-life applications, for instance, large-scale dynamical systems 
and control systems,
the mathematical models are of high complexity and
numerical simulations in high dimensional systems become challenging.
Model order reduction (MOR) is a useful technique in obtaining reasonable 
approximations 
with a significantly reduced computational cost of such problems. 
Through obtaining dominant modes, the dimension of the associated state space is reduced and 
an approximation of the original model with acceptable accuracy is computed.
Proper Orthogonal Decomposition (POD) 
has been used in numerical approximations for dynamic systems 
\cite{volkwein05,Kerschen2005}. 
Recently, POD has been applied to flow problems in porous media 
with high contrasts and heterogeneities 
\cite{globallocal13,global14,poddeim15,cd10,jansen2017use,
trehan2016trajectory,cardoso2009development,efendiev2012local,van2006reduced}. 
The objective of this work is to develop data-driven POD reduced-order models 
for such problems
using  advancements of machine learning. 

A main difficulty of model order reduction is nonlinear problems.
In such problems, modes have to be re-trained or computed in an 
online process \cite{efendiev2016online,nonlinear_AM2015,chung2015residual}, 
which is, in many cases, computationally expensive.
Moreover, in cases with observed data, the 
off-line computed
modes may not honor the data. 
In this work, our goal is to develop an appropriate 
model reduction , which
overcomes the difficulties of nonlinear problems and 
data-present problems. 
We discuss the design of proper orthogonal decomposition (POD) 
modes, where the degrees of freedom have physical
meanings. In particular, they represent the values of the solution
at fixed pre-selected locations. This allows working with POD
coefficients without constructing POD basis functions in cases relevant
to reservoir engineering.

Recently, deep learning has attracted a lot attention (\cite{goodfellow2016deep,lecun2015deep,schmidhuber2015deep}). 
Deep Neural Network (DNN) is at the core of deep learning, and 
is a computing system inspired by the brain's architecture. 
DNN has showed its power in tackling some complex and challenging computer vision tasks, such as classifying a huge number of images, recognizing speech and so on. A deep neural network usually consists of an input layer, an output layer and multiple intermediate hidden layers, with several neurons in each layer. In the learning process, each layer transforms its input data into a little more abstract representation and transfers the signal to the next layer. Some activation functions, acting as the nonlinear transformation on the input signal, are applied to determine whether a neuron is activated or not.

There have been many works discovering the expressivity of deep neural nets theoretically \cite{Cybenko1989, Hornik1991, Csaji2001, Telgrasky2016, Poggio2016, Hanin2017}. The universal approximation property of neural networks has been investigated in a lot of recent studies. It has been shown that deep networks are powerful and versatile in approximating wide classes of functions. Many researchers are inspired to take advantages of the multiple-layer structure of the deep neural networks in approximating complicated functions, and utilize it in the area of solving partial differential equations and model reductions. For instance, in \cite{Ying_paraPDE}, the authors propose a deep neural network to express the physical quantity of interest as a function of random input coefficients, and shows this approach can solve parametric PDE problems accurately and efficiently by some numerical tests. There is another work by E et. al \cite{E_deepRitz}, which aims to represent the trial functions in the Ritz method by deep neural networks. Then the DNN surrogate basis functions are utilized to solve the Poisson problem and eigenvalue problems. In \cite{Shi_resnet}, the authors build a connection between residual networks (ResNet) and the characteristic transport equation. The idea is to propose a continuous flow model for ResNet and show an alternative perspective to understand deep neural networks. 

In this work, we use deep learning concepts combined with 
POD model reduction methodologies constrained at 
observation locations to predict flow dynamics.
We consider a neural network-based approximation of nonlinear flow dynamics.
Flow dynamics is regarded as a multi-layer network, where the solution
at the current time step depends on the solution at the previous time
instant and associated input parameters, such as well rates
and permeability fields.
This allows us to treat the solution via multi-layer network structures,
where each layer is a nonlinear forward map and
 to design novel multi-layer neural network architectures 
for simulations using our reduced-order model concepts.
The resulting forward model takes into account available
data at locations and can be used to reduce the computational
cost associated with forward solves in nonlinear problems.

We will rely on rigorous model reduction concepts to define unknowns 
and connections for each layer. Reduced-order models are important in 
constructing robust learning algorithms since they can identify 
the regions of influence and the appropriate number of variables, thus
allow using small-dimensional maps.
In this work, modified
 proper orthogonal basis functions will be constructed such that 
the degrees of freedom have physical meanings 
(e.g., represent the solution values at selected locations).
Since the constructed basis functions have limited support, 
it will allow localizing the forward dynamics 
by writing the forward map for the solution values at selected locations 
with pre-computed neighborhood structure.
We use a proper orthogonal decomposition model with these 
specifically designed basis functions that are constrained at 
locations.
 A principal component subspace is constructed 
by spanning these basis functions and numerical solutions 
are sought in this subspace. 
As a result,
the neural network is inexpensive to construct.

Our approach combines the available data and physical models,
which constitutes a data-driven modification of the original reduced-order model. 
To be specific, in the network, our reduced-order models will provide a forward map, and 
will also be modified (“trained”) using available observation data. 
Due to the lack of available observation data, we will use computational data to supplement
as needed. The interpolation between data-rich and data-deficient models will also be studied. 
We will also use deep learning algorithms to train the elements of 
the reduced model discrete system. In this case, deep learning architectures 
will be employed to approximate the elements of 
the discrete system and reduced-order model basis functions.

We will present numerical results using deep learning architectures to 
predict the solution and reduced-order model variables. 
In the reduced-order model, designated basis functions allow 
interpolating the solution between observation points. 
A multi-layer neural network based is then built to approximate 
the evolution of the coefficients and, therefore, the flow dynamics. 
We examine how the network architecture, 
which includes the number of layers, and neurons, affects the approximation. 
Our numerical results show that with a fewer number of layers, 
the flow dynamics can be approximated. 
Our numerical results also indicate that the data-driven approach 
improves the quality of approximation.

The paper is organized as follows. In Section \ref{sec:prelim}, we present
a general model and some basic concepts of POD. Section \ref{sec:method}
is devoted to our model learning. In Section \ref{sec:num},
we present numerical results.
We conclude in the last section.

\section{Preliminaries}
\label{sec:prelim}

%

\rev{In this section, we introduce a general problem setting and 
review the concept of POD based global model reduction, 
which is a technique of dimensionality reduction of large-scale 
system of ordinary differential equations (ODE) and its application to 
nonlinear partial differential equations (PDE). 
Consider a time-dependent PDE in the general form 
\begin{equation}
\dfrac{\partial}{\partial t} u = \mathcal{L}(u) + g \quad \text{ in } \Omega \times (0,T), 
\label{eq:pde}
\end{equation}
where $\Omega$ is the spatial domain, $(0,T)$ is the temporal domain, 
$\mathcal{L}$ is a spatio-differential operator on the unknown $u$ 
and $g$ is a given source function. 
The flow dynamic is prescribed to some given initial condition 
and boundary condition. 
We consider spatial discretization procedure by finite element method 
on a Eulerian mesh $\mathcal{T}_h$ for the spatial domain $\Omega$. 
Let $V_h$ be a finite element space spanned by the nodal basis $\{ \phi_j \}_{j=1}^n$ on $\mathcal{T}_h$. 
We seek numerical solution of \eqref{eq:pde} by an expansion 
\begin{equation}
u(x,t) = \sum_{j=1}^n y_j(t) \phi_j(x), 
\end{equation}
which yields a system of ODE in the form 
\begin{equation}
\dfrac{d}{dt} \mathbf{y}(t) = \mathbf{B} \mathbf{y}(t)  + \mathbf{f}(\mathbf{y}(t)),
\end{equation}
where $\mathbf{y}(t) \in \mathbb{R}^n$ is the state vector, 
$\mathbf{B} \in \mathbb{R}^{n \times n}$ is a constant matrix, 
and $\mathbf{f}: \mathbb{R}^n \to \mathbb{R}^n$ is a nonlinear function. 
In our applications, the dimension $n$  
corresponds to the number of physical grid points in the mesh. 
In general, the dimension $n$ is huge and model reduction techniques 
provide efficient reduced-order models and bring computational savings.  
}

\subsection{Proper Orthogonal Decomposition}

Proper Orthogonal Decomposition is a popular mode decomposition method, 
which aims at reducing the order of the model by 
extracting important relevant feature representation with a 
low dimensional space.
In this section, we briefly discuss the POD method. 
For a more detailed discussion of the use of POD on dynamic systems, 
the reader is referred to \cite{volkwein05,Kerschen2005}.
In POD, a low-dimensional set of modes, i.e., important degrees of freedom, 
are identified based on processing information 
from a sequence of snapshots, 
i.e., instantaneous solutions from the dynamic process,
and extracting the most energetic structures in terms of the largest singular values. 
In the statistical point of view, the extracted modes are uncorrelated and 
form an optimal reduced order model, in the sense that the variance is maximized and 
the mean squared distance between the snapshots and the POD subspace is minimized. 

\rev{
Proper orthogonal decomposition starts with a collection of $N \ll n$ instantaneous snapshots 
$\{ \mathbf{y}_j \}_{j=1}^N \subset \mathbb{R}^n$, 
where the snapshot times in the above sequence is assumed to be equidistant. 
The snapshots span a snapshot space of dimension $r$ 
and are arranged in a matrix form known as the snapshot matrix  
\begin{equation}
\mathbf{Y} = [\mathbf{y}_1 \;\; \mathbf{y}_2 \;\; \cdots \;\; \mathbf{y}_N] \in \mathbb{R}^{n \times N}. 
\end{equation}
The idea of POD is to seek the subspace of a certain dimension 
which best approximates the linear space spanned by the snapshots. 
Among all subsets of $m < r$ orthonormal vectors in $\mathbb{R}^n$, 
we seek the POD basis $\{ \mathbf{v}_j\}_{j=1}^m$ by solving a minimization problem 
\begin{equation}
\argmin_{ \substack{ \{ \mathbf{v}_j\}_{j=1}^m \subset \mathbb{R}^n \\ \langle \mathbf{v}_i, \mathbf{v}_j \rangle = \delta_{ij} }  }
\sum_{i=1}^N \left\| \mathbf{y}_i - \sum_{j=1}^m \left\langle \mathbf{y}_i, \mathbf{v}_j \right\rangle \mathbf{v}_j \right\|_2^2, 
\label{eq:podmin}
\end{equation}
The minimzation problem is processed by performing a singular value decomposition 
on the snapshot matrix $\mathbf{Y}$ 
\begin{equation}
\mathbf{Y} = \mathbf{V} \Lambda \mathbf{W}^T, 
\end{equation}
where $\mathbf{V} = [\mathbf{v}_1, \mathbf{v}_2, \cdots, \mathbf{v}_r] \in \mathbb{R}^{n \times r}$ 
and $\mathbf{W} = [\mathbf{w}_1, \mathbf{w}_2, \cdots, \mathbf{w}_r] \in \mathbb{R}^{N \times r}$
consist of the left-singular vectors and right-singular vectors of $\mathbf{Y}$ respectively,
and $\Lambda = \text{diag}(\sigma_1, \sigma_2, \cdots, \sigma_r) \in \mathbb{R}^{r \times r}$ 
is the diagonal matrix consisting of the singular values of $\mathbf{Y}$. 
Constructively, we denote the correlation matrix from the snapshot sequence by $\mathbf{C} = \mathbf{Y}^T \mathbf{Y}$, 
and compute the eigenvalue decomposition on $\mathbf{C}$
\begin{equation}
\mathbf{C}\mathbf{q}_j = \lambda_j \mathbf{q}_j,
\end{equation}
and obtain the singular values $\{ \sigma_j \}_{j=1}^r$ and singular vectors $\{\mathbf{v}_j\}_{j=1}^r$ by
\begin{equation}
\sigma_j = \sqrt{\lambda_j} \text{ and } \mathbf{v}_j = \dfrac{1}{\sigma_j} \mathbf{Y} \mathbf{q}_j.
\end{equation}
Here the singular values are arranged in descending energy ranking, i.e., 
$\sigma_1 \geq \sigma_2 \geq \cdots \geq \sigma_r > 0$, 
which correspond to the energy content of a mode. 
The energy ranking provides a measure of the importance of the mode 
in capturing the relevant dynamic process. 
The POD basis, i.e. the solution of the minimization problem \eqref{eq:podmin}, 
is then given by selecting the first $m$ singular vectors $\{\mathbf{v}_j\}_{j=1}^m$. 
In this case, we have 
\begin{equation}
\sum_{i=1}^N \left\| \mathbf{y}_i - \sum_{j=1}^m \left\langle \mathbf{y}_i, \mathbf{v}_j \right\rangle \mathbf{v}_j \right\|_2^2 
= \sum_{j=m+1}^r \sigma_j^2. 
\end{equation}
The size $m$ of the POD basis has to be sufficiently large to 
include the first few largest singular values and 
ensure a good approximation to the snapshot matrix. 
The number of basis can be pre-defined or determined by means of fractional energy, i.e. 
fixing a threshold $E_0$, pick the smallest integer $m$ such that 
\begin{equation}
E = \dfrac{\sum_{j=1}^m \sigma_j^2}{\sum_{j=1}^r \sigma_j^2} > E_0,
\end{equation}
In general, a few basis is needed if the singular values decay quickly. 
The rate of decay depends on the intrinsic dynamics of the system and the selection of the snapshots.
}  

\rev{
\subsection{Fully discrete reduced-order model}
Using the aforementioned POD basis $\{ \mathbf{v}_j \}_{j=1}^m$, 
we can express the solution as
\begin{equation}
\mathbf{y}(t) \approx \sum_{j=1}^m \widetilde{c}_j(t) \mathbf{v}_j = \mathbf{V} \widetilde{c}(t),
\end{equation}
where $\widetilde{c}(t) = \left( \widetilde{c}_1(t), \widetilde{c}_2(t), \ldots, \widetilde{c}_m(t) \right)^T \in \mathbb{R}^m$ 
is the coordinates of $\mathbf{y}(t)$ with respect to the POD basis. 
We therefore derive a reduced-order ODE system 
\begin{equation}
\dfrac{d}{dt} \widetilde{c}(t) = \mathbf{V}^T \mathbf{B} \mathbf{V} \widetilde{c}(t)  + \mathbf{V}^T \mathbf{f}(\mathbf{V} \widetilde{c}(t)), 
\end{equation}
and further reduce it to an algebraic system. 
We consider a partition $0 = t_0 < t_1 < \ldots < t_s = T$ for the temporal domain $(0,T)$. 
Using, for example, implicit Euler method for temporal discretization, we obtain a recurrence relation 
\begin{equation}
\widetilde{c}^{n+1} = \widetilde{c}^n + 
(t_{n+1} - t_n) \left( \mathbf{V}^T \mathbf{B} \mathbf{V} \widetilde{c}^{n+1}  + \mathbf{V}^T \mathbf{f}(\mathbf{V}\widetilde{c}^{n+1}) \right), 
\end{equation}
where $\widetilde{c}^{n}$ denotes the numerical solution of $\widetilde{c}(t)$ at the time instant $t = t_n$. 
The nonlinear term $\mathbf{f}$ can be handled with different techniques, such as 
direct linearization method, fixed point iterations and Discrete Empirical Interpolation Method (DEIM), 
depending on situations and need for accuracy in particular applications. 
}

\subsection{Construction of nodal basis functions}

Next, we present the construction of basis functions \rev{in the POD subspace}. 
The basis functions are designed such that the degrees
of freedom have physical meanings 
(e.g., represent the solution values at selected locations).
Since the constructed basis functions have limited support, it will allow 
localizing the forward dynamics 
by writing the forward map for the solution values at selected locations 
with pre-computed neighborhood structure.

Given a set of nodes $\{ x_k \}_{k=1}^m$ in the mesh $\mathcal{T}_h$, 
which correspond to particular physical points in the spatial domain $\Omega$, 
we construct nodal basis functions 
by linear combinations of POD modes $\{\mathbf{v}_j\}_{j=1}^m$. 
More precisely, we seek coefficients $\alpha_{ij}$ such that 
\begin{equation}
\sum_{j=1}^m \alpha_{ij} \rev{\mathbf{V}_{kj}} = \delta_{ik}.
\end{equation}
\rev{We remark that $\mathbf{V}_{kj}$ is the nodal evaluation 
of the interpolant of $\mathbf{v}_j$ 
in the finite element space $V_h$ at the node $x_k$.}
The nodal basis $\psi_k$ \rev{at the node $x_k$} is then defined by 
\begin{equation}
\psi_k = \sum_{j=1}^m \alpha_{kj} \mathbf{v}_j.
\end{equation}
\rev{The set of nodal basis spans exactly 
the POD subspace. 
A reduced-order state vector $\mathbf{y}(t)$ in the POD subspace 
can then be written in the expansion 
\begin{equation}
\mathbf{y}(t) = \sum_{k=1}^m c_k(t) \psi_k,
\end{equation}
where the coefficients $c_k(t)$ represent the nodal evaluation 
of the finite element approximation of $u(x,t)$ at the node $x_k$. 
Furthermore, the coefficients $\widetilde{c}^n$ of the original POD basis functions and $c^n$ of 
the POD nodal basis functions are related by 
\begin{equation}
c^n = \mathbf{V}_m \widetilde{c}^n,
\end{equation}
where $\mathbf{V}_m \in \mathbb{R}^{m \times m}$ is 
the submatrix obtained from $\mathbf{V}$ by 
taking the rows corresponding to the nodes $\{ x_k \}_{k=1}^m$.}

\section{Deep Global Model Reduction and Learning}
\label{sec:method}

\subsection{Main idea}
We will make use of the reduced-order model described in Section~\ref{sec:prelim} to model the flow dynamics, 
and a deep neural network to approximate the flow profile. 
In many cases, the flow profile is dependent on data. 
The idea of this work is to make use of deep learning to 
combine the reduced-order model and 
available data and provide an efficient numerical model 
for modelling the flow profile.

First, we note that the solution at the time instant $n+1$ depends on the solution at the time instant $n$ 
and input parameters $I^{n+1}$, such as permeability field and source terms. 
Here, we would like to use a neural network to describe the relationship of the solutions between two consecutive time instants. Suppose we have a total $m$ sample realization in the training set. 
For each realization, given a set of input parameters, 
we solve the aforementioned reduced-order model 
and obtain the coefficients at particular points
\begin{equation} \{c^{0}, \cdots,  c^{k}\}
\end{equation}
at all time steps. Our goal is to use deep learning techniques to train the trajectories  
and find a network $\mathcal{N}$ to describe the pushforward map 
between $c^{n}$ and $c^{n+1}$ for any training sample.
\begin{equation}
c^{n+1} 	\sim \mathcal{N} (c^n, I^{n+1}),
\end{equation}
where $I^{n+1}$ is an input parameter which could vary over time, 
and $\mathcal{N}$ is a multi-layer network to be trained.
The network $\mathcal{N}$ will approximate the discrete flow dynamics. 

In our neural network, $c^n$ and $I^{n+1}$ are the inputs, $c^{n+1}$ is the output. 
One can take the coefficients from time $0$ to time $k-1$ as input, 
and from time $1$ to $k$ as output in the training process. 
In this case, a universal neural net $\mathcal{N}$ is obtained. 
The solution at time $0$ can then be forwarded all the way to time $k$ 
by repeatedly applying the universal network $k$ times, that is,
\begin{equation}
c^{k} 	\sim  \mathcal{N} (\mathcal{N} \cdots \mathcal{N} (c^0, I^1) \cdots , I^{k-1}), I^{k}).
\end{equation}
After a network is trained, it can be used for predicting the trajectory given a new set of input parameters $I^{n+1}$ and 
realization of coefficients at initial time 
$c_{\text{new}}^0$ by
\begin{equation}
c_{\text{new}}^{k} 	\sim  \mathcal{N} (\mathcal{N} \cdots \mathcal{N} (c_{\text{new}}^0, I^1) \cdots , I^{k-1}), I^k).
\end{equation}
Alternatively, one can also train each forward map for any two consecutive time instants as needed. That is, we will have $c^{n+1} \sim  \mathcal{N}_{n+1}(c^n , I^{n+1})$, for $n =0, 1,\ldots, k-1 $.  In this case, to predict the final time solution ${c_{\text{new}}^{k}}$ given the initial time solution ${c_{\text{new}}}^{0}$, we use $k$ different networks $\mathcal{N}_1, \cdots, \mathcal{N}_k$
\begin{equation}
c_{\text{new}}^{k} 	\sim  \mathcal{N}_k (\mathcal{N}_{k-1} \cdots \mathcal{N}_1 (c_{\text{new}}^0, I^1) \cdots , I^{k-1}), I^k).
\end{equation}
We remark that, besides the solution $u^n$ at the previous time instant, the other input parameters $I^{n+1}$ such as permeability or source terms can be different when entering the network at different time steps.

In this work, we would like to incorporate available observed data in the neural network. The observation data will help to supplement the computational data which are obtained from the underlying reduced order model, and improve the performance of the neural network model such that it will take into account real data effects. From now on, we use
 $\{c_{s}^{0},\cdots,  c_{s}^{k}\}$ to denote the simulation data, and $\{c_{o}^{0},\cdots,  c_{o}^{k}\}$ to denote the observation data. 
 
One can get the observation data from real field experiment. However in this work, we generate the observation data by running a new simulation on the ``true permeability field'' using standard finite element method, and using the results as observed data. For the computational data, we will perturb the ``true permeability field'', and use the reduced-order model, i.e., POD model for simulation. In the training process, we are interested in investigating the effects of observation data in the output. One can compare the performance of deep neural networks when using different combinations of computation and observation data.

For the comparison, we will consider the following three networks
\begin{itemize}
\item Network A: Use all observation data as output,
\begin{equation}\label{eq:No}
{c_o}^{n+1} \sim \mathcal{N}_{o} ({c_s}^n, I^{n+1})
\end{equation}

\item Network B: Use a mixture of observation data and simulation data as output,
\begin{equation}\label{eq:Nm}
{c_m}^{n+1} \sim \mathcal{N}_m ({c_s}^n, I^{n+1})
\end{equation}

\item Network C: Use all simulation data (no observation data) as output,
\begin{equation}\label{eq:Ns}
{c_s}^{n+1} \sim \mathcal{N}_s ({c_s}^n, I^{n+1})
\end{equation}

\end{itemize}
where $c_m$ is a mixture of simulation data and observed data.

The first network (Network A) corresponds to the case when the observation data is sufficient. One can merely utilize the observation data in the training process. That is, the observation data at time $n+1$ can be learnt as a function of the observation data at time $n$. This map will fit the real data very well given enough training data; however, it will not be able to approximate the reduced-order model. Moreover, in the real application, the observation data are hard to obtain, and in order to make the training effective, deep learning requires a huge amount of data. Thus Network A is not applicable in real case, and we will use the results from Network A as a reference.

The third network (Network C), on the other hand, will simply take all simulation data in the training process. In this case, one will get a network describes the simulation model (in our example, the POD reduced-order model) as best as it can but ignore the observational data effects. This network can serve as an emulator to do a fast simulation. We will also utilize Network C results as a reference.

 We are interested in investigating the performance of Network B, where we take a combination of computational data and observational data to train. It will not only take in to account the underlying physics but also use the real data to modify the reduced-order model, thus resulting in a data-driven model.

\subsection{Network structures}

\begin{figure}[ht!]
\centering
\includegraphics[width=0.48\linewidth]{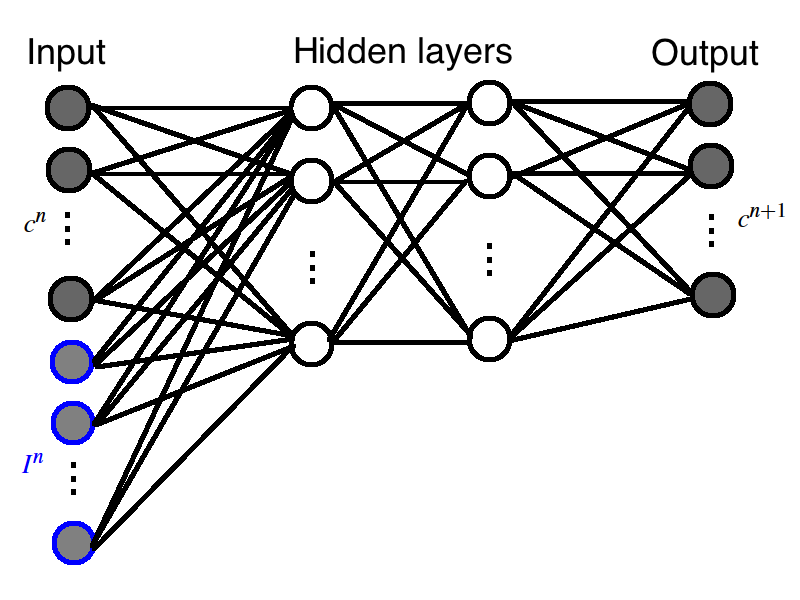}
\caption{An illustration of deep neural network}
\label{fig:dnn}
\end{figure}

Mathematically, a neural network $\mathcal{N}$ of $L$ layers with 
input $\mathbf{x}$ and output $\mathbf{y}$ 
is a function in the form
\begin{equation*}
\mathcal{N}(\mathbf{x}; \theta) = \sigma(W_L \sigma (\cdots  \sigma(W_2  \sigma(W_1 \mathbf{x} + b_1) + b_2) \cdots   ) + b_L),
\end{equation*} 
where $\theta : = (W_1, \cdots, W_L, b_1, \cdots, b_L)$ is a set of network parameters, 
$W$'s are the weight matrices and $b$'s are the bias vectors. The activation function $\sigma$ acts as entry-wise evaluation. 
A neural network describes the connection of a collection of nodes (neurons) sit in successive layers. 
The output neurons in each layer is simultaneously the input neurons in the next layer. 
The data propagate from the input layer to the output layer through hidden layers. 
The neurons can be switched on or off as the input is propagated forward through the network.
The weight matrices $W$'s control the connectivity of the neurons. 
The number of layers $L$ describes the depth of the neural network.
Figure~\ref{fig:dnn} depicts a deep neural network in out setting, 
in which each circular node represents a neuron 
and each line represents a connection from one neuron to another. 
The input layer of the neural network consists of 
the coefficients $c^n$ and the input parameters $I^n$. 

Given a set of data $(\mathbf{x}_j, \mathbf{y}_j)$, the deep neural network aims to find the parameters $\theta^*$ by solving an optimization problem
\begin{equation*}
\theta^* = \argmin_{\theta} \frac{1}{N}\sum_{j=1}^{N} ||\mathbf{y}_j - \mathcal{N}(\mathbf{x}_j; \theta) ||^2_2,
\end{equation*}
where $N$ is the number of the samples. Here, the function $L(\theta) = \frac{1}{N}\sum_{j=1}^{N} ||\mathbf{y}_j - \mathcal{N}(\mathbf{x}_j; \theta) ||^2_2$ is known as the loss function. One needs to select suitable number of layers, number of neurons in each layer, the activation function, the loss function and the optimizers for the network. 

As discussed in the previous section, we consider three different networks, 
namely $\mathcal{N}_o$, $\mathcal{N}_{m}$ and $\mathcal{N}_s$. 
For each of these networks, we take the vector $\mathbf{x} = ({c_s}^n, I^{n+1})$ 
containing the numerical solution vectors and the data at a particular time step as the input. 
In our setting, the input parameter $I^{n+1}$, if present, could be the static permeability field or 
the source function.
Based on the availability of the observational data in the sample pairs, 
we will select an appropriate network among \eqref{eq:No}, \eqref{eq:Nm} and \eqref{eq:Ns} accordingly.
The output $\mathbf{y} = {c_\alpha}^{n+1}$ is taken as the numerical 
solution at the next time instant,
where $\alpha = o,m,s$ corresponds to the network. 

Here, we briefly summarize the architecture of 
the network $\mathcal{N}_\alpha$, 
where $\alpha = o,m,s$ for three networks we defined in \eqref{eq:No}, \eqref{eq:Nm} and \eqref{eq:Ns} respectively. 

As for the input of the network, we use $\mathbf{x} = ({c_s}^n, I^{n+1})$, which are the vectors containing the numerical solution vectors and the input parameters in a particular time step. The corresponding output data are $\mathbf{y} = {c_\alpha}^{n+1}$, which contains the numerical solution in the next time step. In between the input and output layer, we test on 3--10 hidden layers with 20-400 neurons in each hidden layer. In the training, there are $N = mk$ sample pairs of $(\mathbf{x}_j, \mathbf{y}_j)$ collected, where $m$ is the number of realizations of flow dynamics and $k$ is the number of time steps.

In between layers, we need the activation function. The ReLU function (rectified linear unit activation function) is 
a popular choice for activation function in training deep neural network architectures \cite{glorot11}. 
However, in optimizing a neural network with ReLU as activation function, 
weights on neurons which do not activate initially will not be adjusted, 
resulting in slow convergence.  Alternatively, leaky ReLU can be employed to avoid such scenarios \cite{relu}. 
We choose leaky ReLU in our network structure. 
As for the training optimizer, we use AdaMax \cite{adam}, which is a stochastic gradient descent (SGD) type algorithm well-suited for high-dimensional parameter space, in minimizing the loss function.

\begin{remark}
We note that the discretization parameters can be learned using proposed DL techniques. Our numerical results show that one can achieve a very good accuracy (about 2 \%) using DNN with 5 layers and 100 neurons in each layer.
These results show that deep learning techniques can be used for fast prediction of discretized systems. In our paper, our main goal was to incorporate the observed data into the dynamical system and modify reduced-order models to take
into account the observed data.
We have also used deep learning techniques to predict online basis functions \cite{efendiev2016online,nonlinear_AM2015,chung2015residual}, which allow re-constructing the fine-scale solution.
\end{remark}

\section{Numerical examples}
\label{sec:num}

\rev{
In this section, we present numerical examples. 
We apply our method to predict the evolution 
of the pressure in a nonlinear single-phase flow problem, 
and the saturation in a a two-phase flow problem. 
Using POD global model reduction technique, 
we obtain coefficients of numerical solutions in the reduced-order model 
and use as training samples to construct neural network approximations 
of the corresponding nonlinear flow dynamics. 
All the network training are performed using the 
Python deep learning API Keras \cite{chollet2015keras}.
}

\subsection{Application to single-phase flow} 

\begin{figure}[ht!]
\centering
\includegraphics[width=0.3\linewidth]{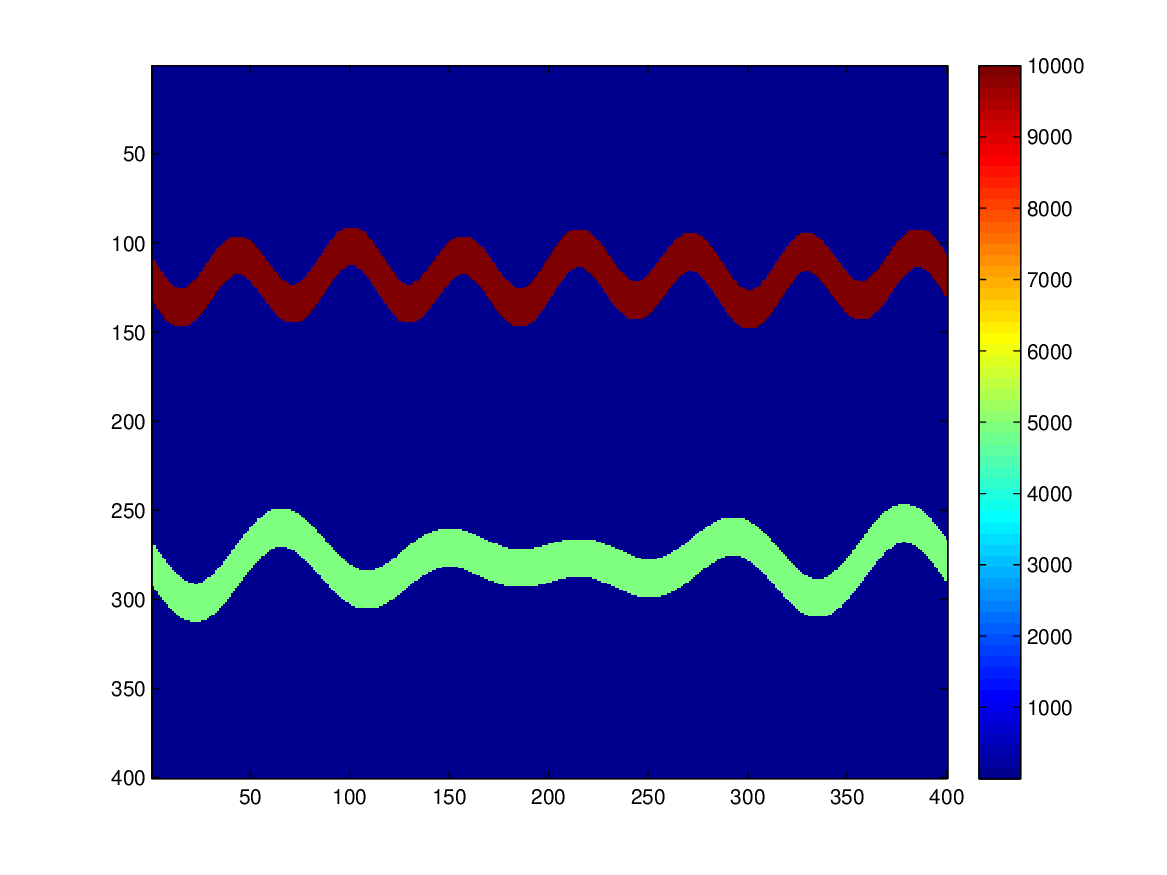}
\includegraphics[width=0.3\linewidth]{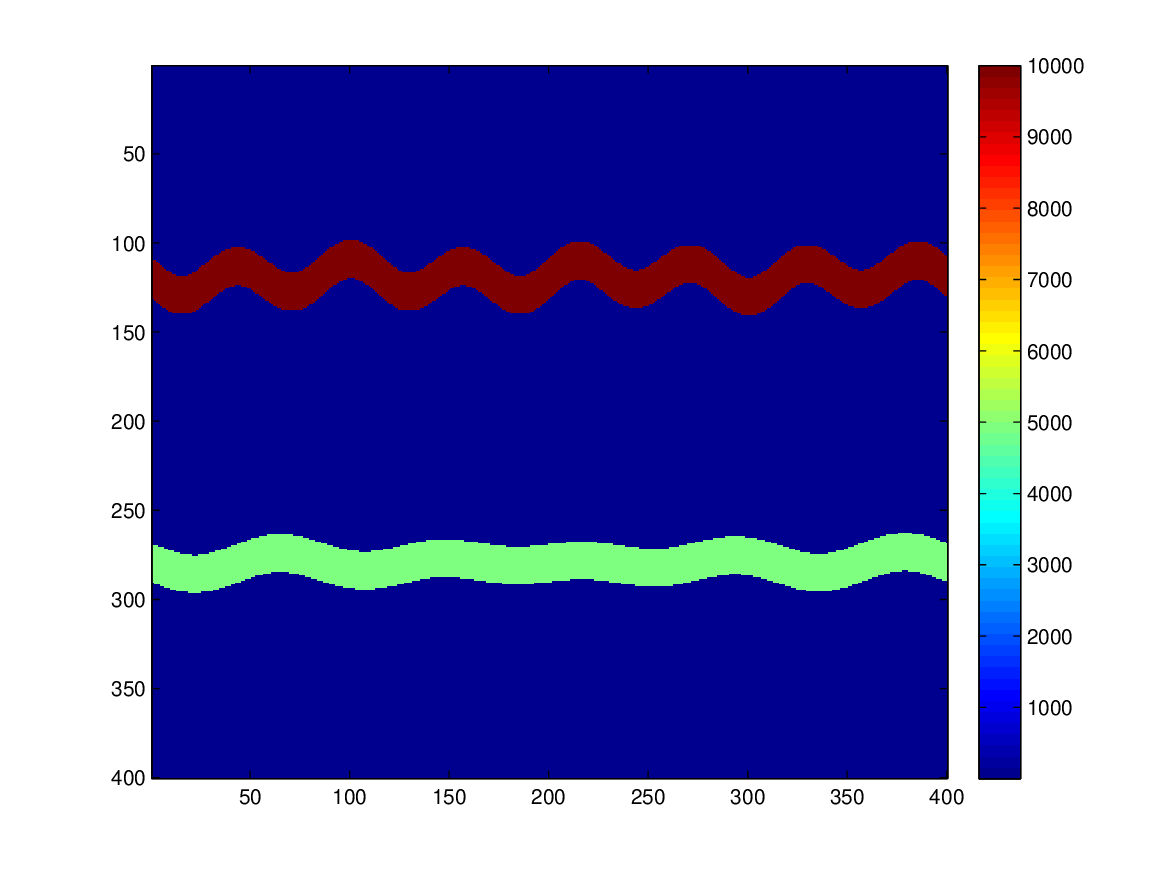}
\includegraphics[width=0.3\linewidth]{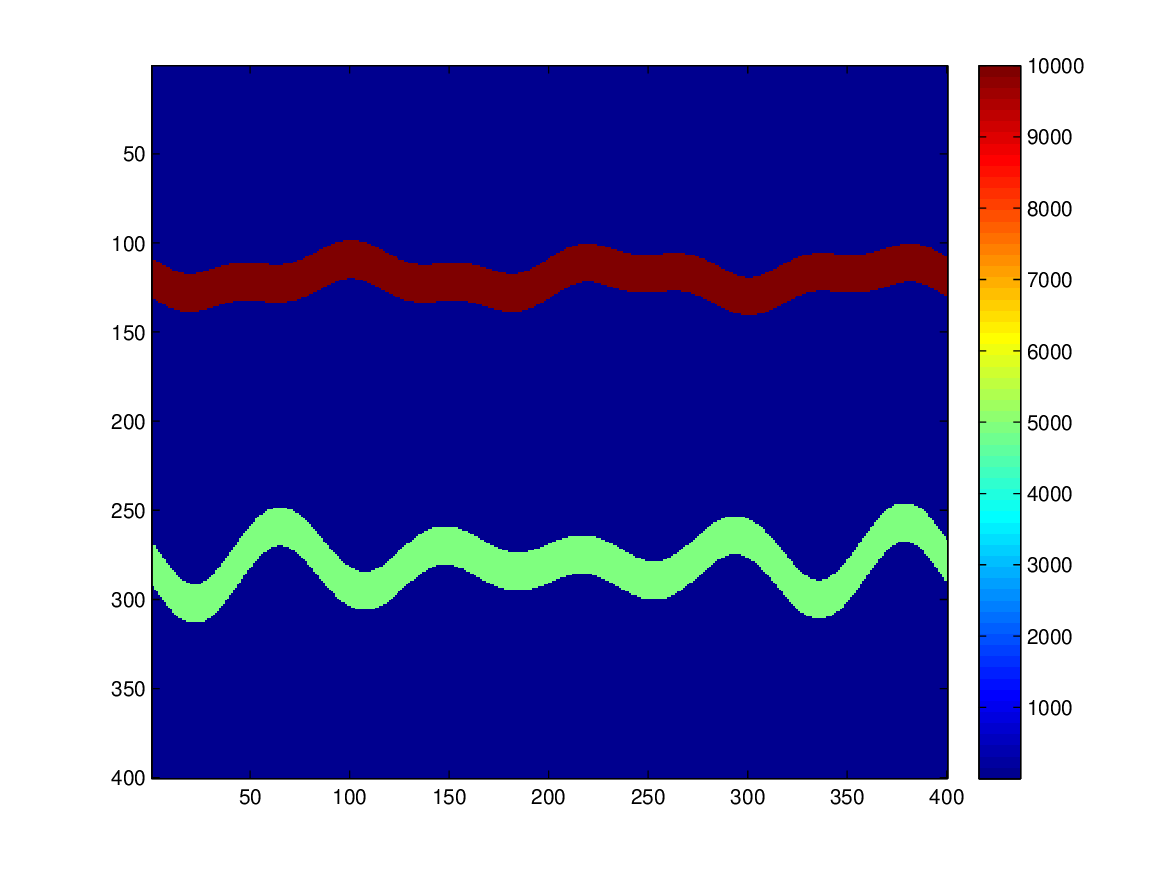}
\caption{Samples of static permeability field used in single-phase flow.}
\label{fig:kappa}
\end{figure}

As a first example, we consider a simple nonlinear single-phase flow 
in the spatial domain $\Omega = [0,1] \times [0,1]$:
\begin{equation}\label{eq:diffusion}
\frac{\partial u}{\partial t}- \text{div} (\kappa(x,u)  \nabla u) = g \quad \text{in} \quad \Omega,
\end{equation}
subject to homogeneous Dirichlet boundary 
condition $u\vert_{\partial \Omega} = 0$. This equation describes
unsaturated flow in heterogeneous media,
which are widely used \cite{richards1931capillary,gardner1958some,van1980closed,dostert2009efficient,celia1990general}.
In our simulations, we will use an exponential model
$\kappa(x,u)=\kappa(x)\exp(\alpha u)$.
Here, $u$ is the pressure of flow, $g$ is a time-dependent source term 
and $\alpha$ is a nonlinearity parameter.
The function $\kappa(x)$ is a stationary heterogeneous permeability field 
of high contrast, i.e., with large variations within the domain $\Omega$. 
In this example, we focus on permeability fields that contain wavelet-like 
channels as shown in Figure~\ref{fig:kappa}. In each realization of 
the permeability field, 
there are two non-overlapping channels with high conductivity values in the domain $\Omega$, 
while the conductivity value in the background is $1$. 
Channelized permeability fields
are challenging for model reduction and prediction and, thus, we focus
on flows corresponding to these permeability fields. The numerical
tests for Gaussian permeability fields
(\cite{efendiev2005efficient,vo2014new})
show a good accuracy because of the smoothness of the solution with respect
to the parameters.

Next, we present the details of numerical discretization of the problem. 
Suppose the spatial domain $\Omega$ is partitioned into a rectangular mesh $\mathcal{T}_h$, 
and a set of piecewise bilinear conforming finite element basis functions $\{ v_j \}$
is constructed on the mesh. 
We denote the finite element space by $V_h = \mathbb{Q}^1(\mathcal{T}_h)$.
Using direct linearization for the nonlinear term, 
implicit Euler method for temporal discretization and 
a Galerkin finite element method for spatial discretization, 
the numerical solution $u_h^{n+1}$ at the time instant $n+1$ 
is obtained by solving the following variational formulation: 
find $u_h^{n+1} \in V_h$ such that
\begin{equation}\label{eq:fine}
\int_\Omega \frac{u_h^{n+1} - u_h^n }{\Delta t} v  +  \int_\Omega \kappa \exp(\alpha u_h^{n}) \nabla u_h^{n+1} \cdot \nabla v = \int_{\Omega} g^{n+1} v \quad \text{ for all } v \in V_h.
\end{equation}
Here $\Delta t$ is the time step and $h$ is the mesh size.
With a slight abuse of notation, we again denote the coefficients of 
the numerical solution with the piecewise bilinear basis functions by $u_h^{n+1}$. 
Then, the variational formulation can be written in the matrix form
\begin{equation}
u_h^{n+1} =  (M + \Delta t A(u_h^{n}))^{-1} (M u_h^{n} + \Delta t b^{n+1}),
\end{equation}
where $M$, $A(u_h^{n})$ and $b^{n+1}$ 
are the mass matrix, the stiffness matrix and the load vector 
with respect to the bilinear basis functions $v_j$, i.e., 
\begin{equation}
\begin{split}
M_{ij} & = \int_\Omega v_i v_j, \\
[A(u_h^{n})]_{ij} & = \int_\Omega \kappa \exp(\alpha u_h^{n}) \nabla v_i \cdot \nabla v_j, \\
b^{n+1}_i & = \int_\Omega g^{n+1} v_i. 
\end{split}
\end{equation}

\begin{figure}[ht!]
\centering
\includegraphics[width=0.3\linewidth]{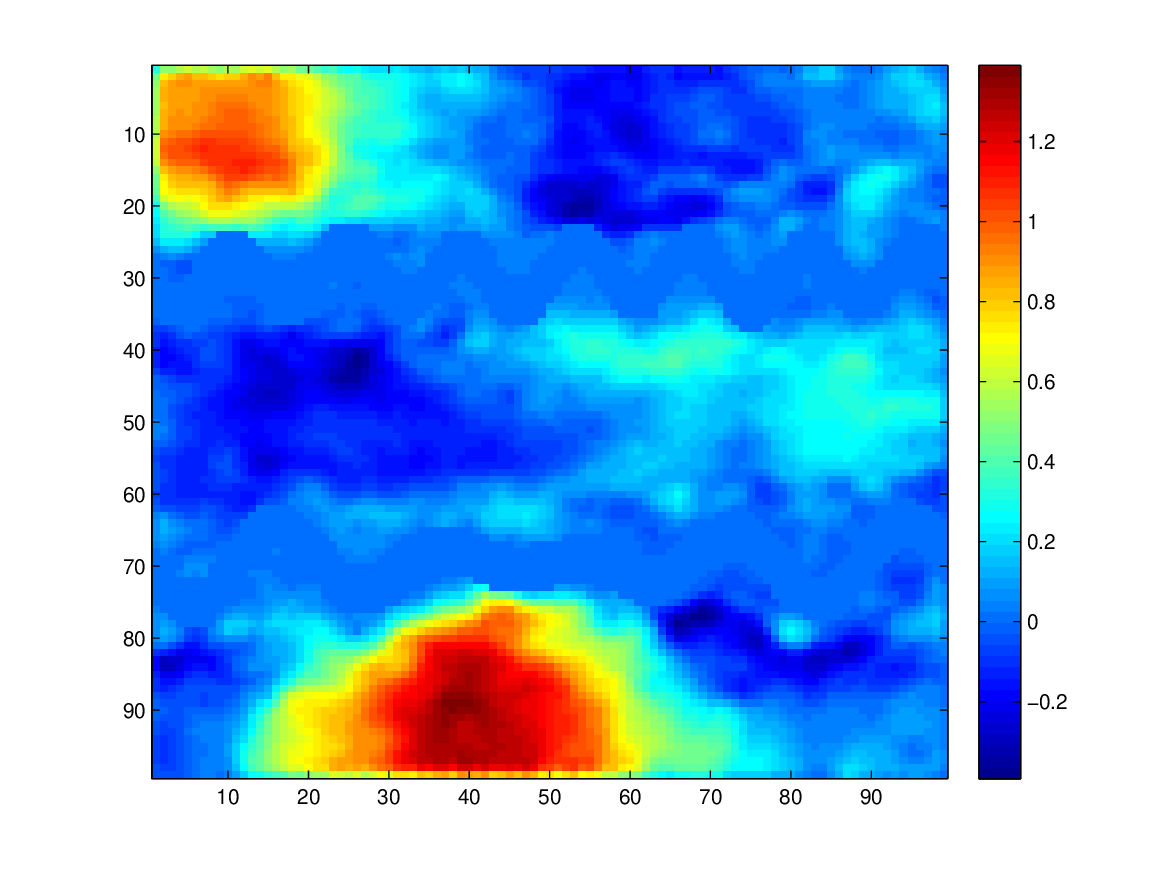}
\includegraphics[width=0.3\linewidth]{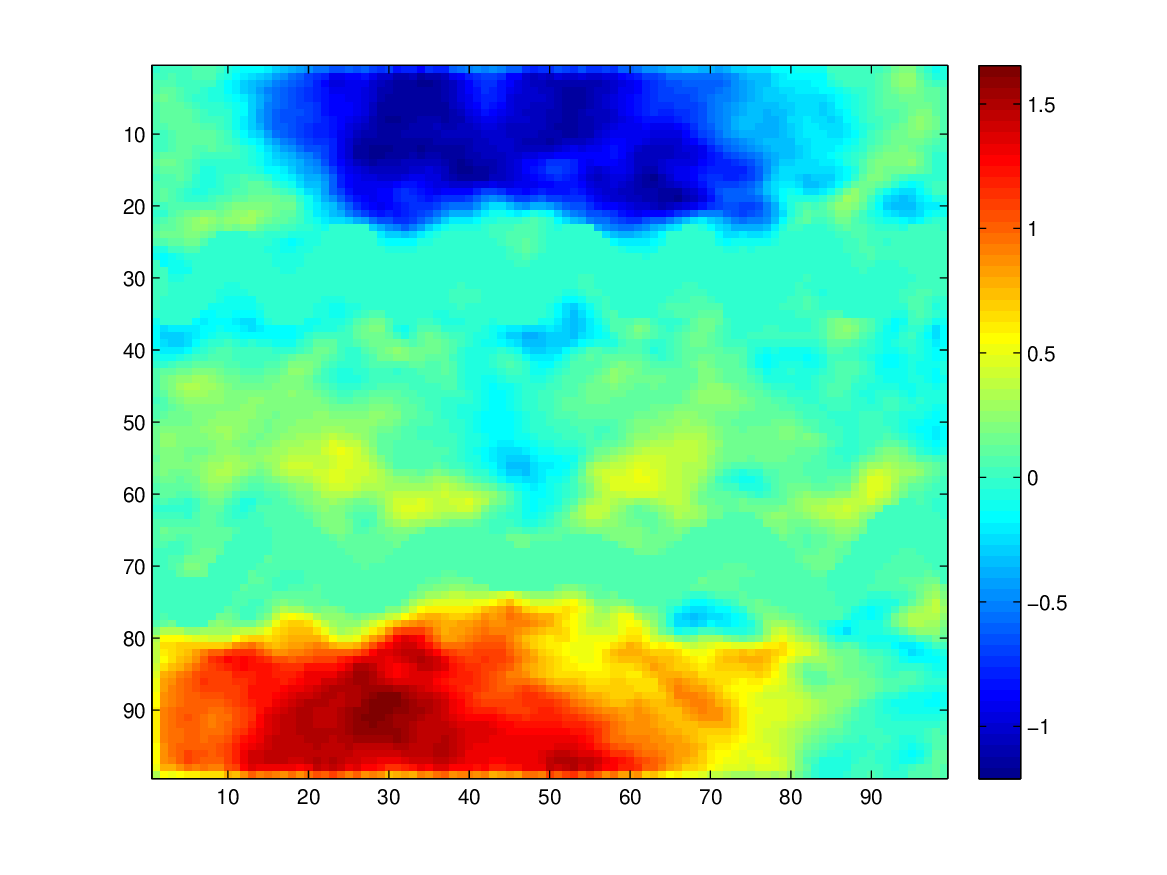}
\includegraphics[width=0.3\linewidth]{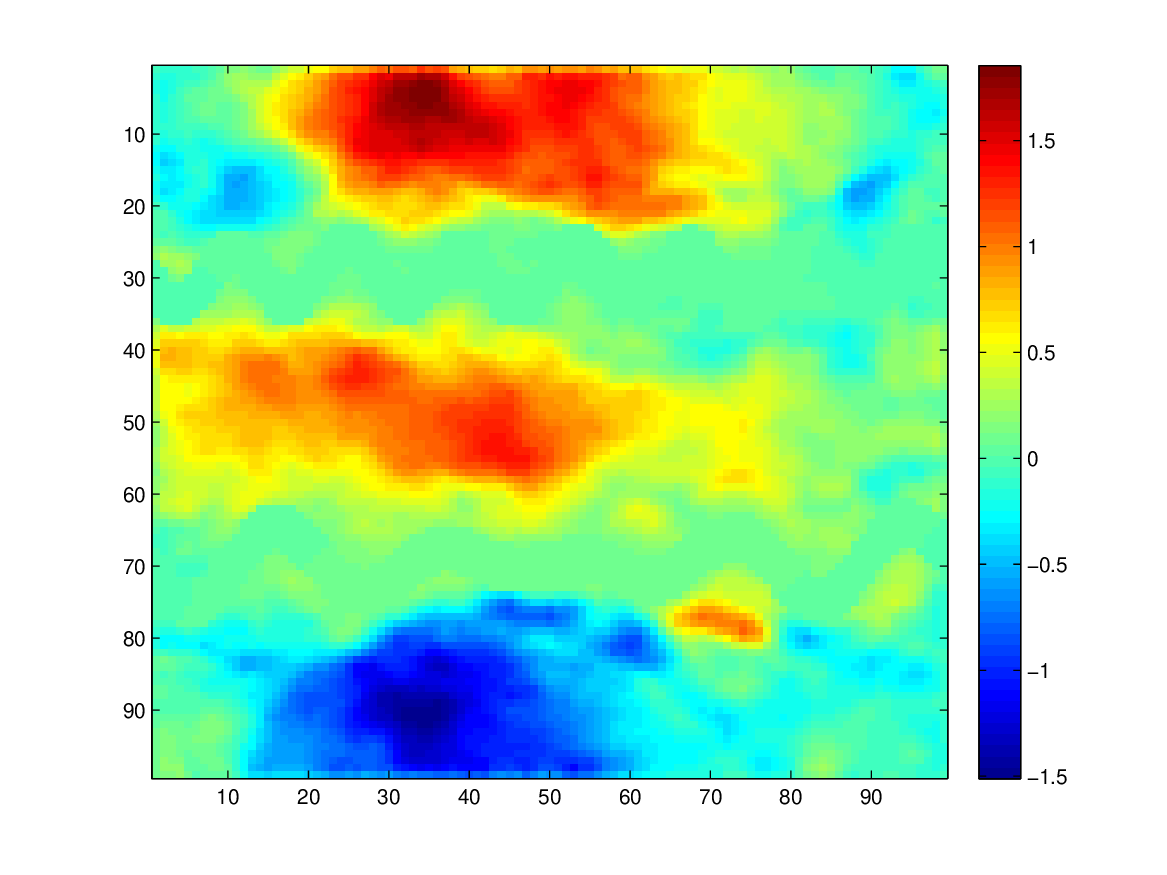}
\newline
\includegraphics[width=0.3\linewidth]{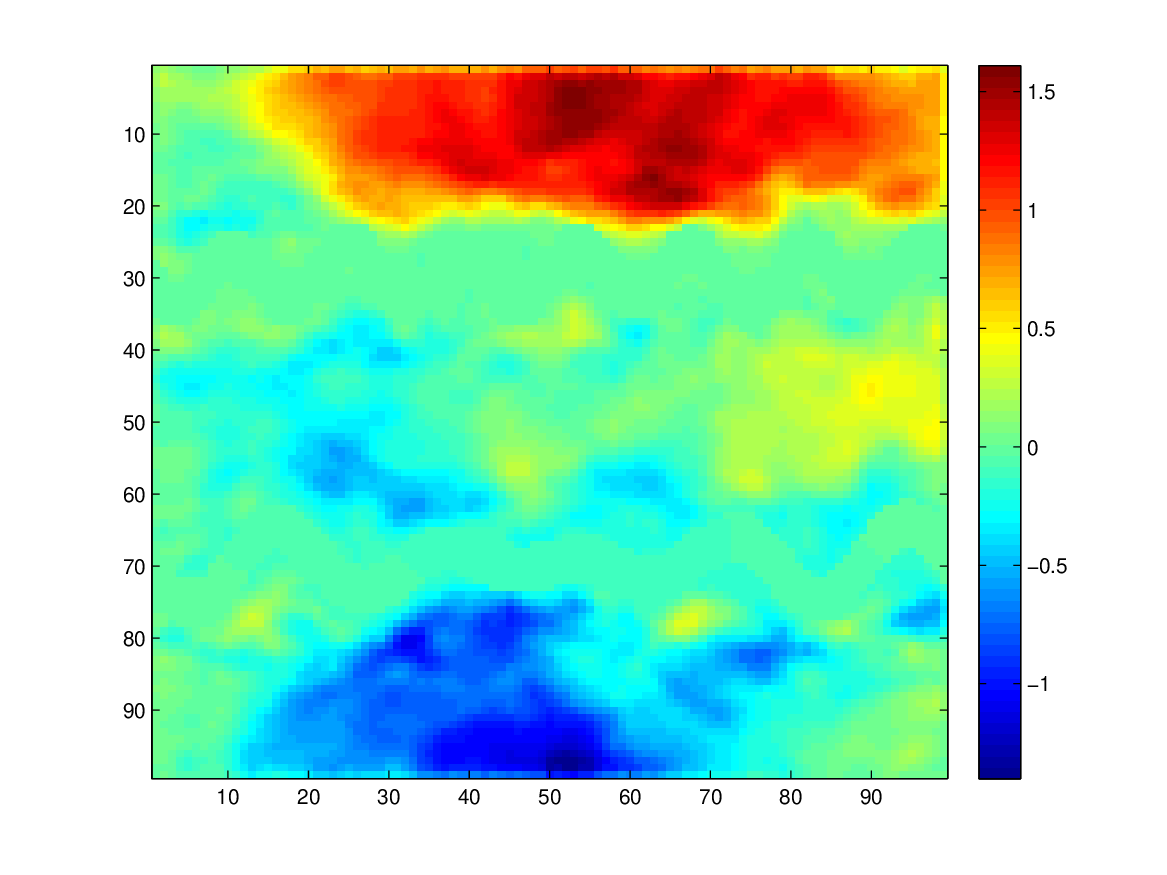}
\includegraphics[width=0.3\linewidth]{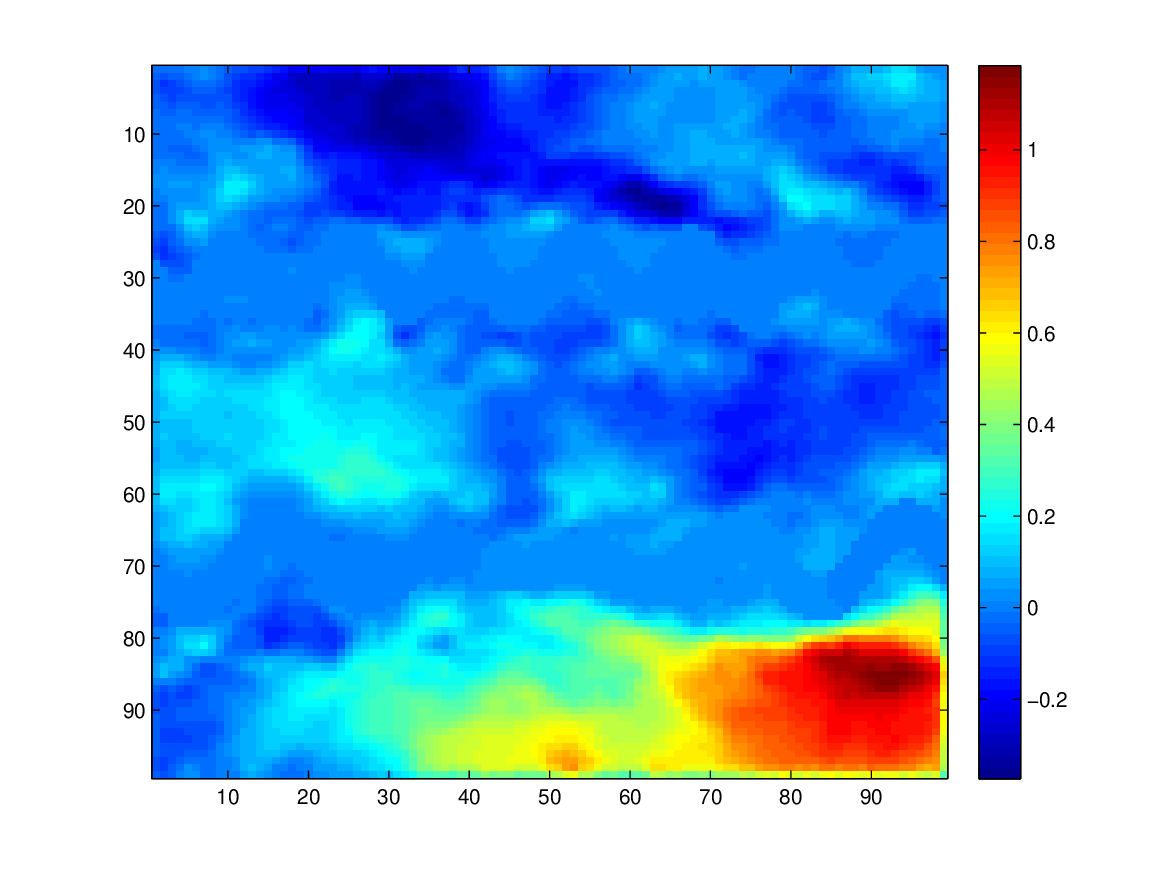}
\caption{Illustration of nodal basis functions.}
\label{fig:nodal}
\end{figure}

In our simulation, the flow is simulated from an initial time $t = 0$ 
to a final time $t = 0.01$ in $10$ time steps. 
Realizations of flow dynamics are computed using 
independent and uniformly distributed initial conditions. 
We use POD to extract dominant modes from snapshot solutions 
and construct POD nodal basis functions. 
Examples of POD nodal basis functions are shown in Figure~\ref{fig:nodal}. 
Simulation data of the dynamic process under the reduced-order model 
are then obtained and used in the training set. 
Different forms of inputs, depending on situations, are investigated. 
Using these data as samples, universal multi-layer networks are trained to 
approximate the flow dynamics. 
We use the trained networks to predict the output with some 
new unseen inputs, and reconstruct the numerical solution 
using the predicted coefficients. 
We examine the quality of our networks by computing the $L^2$ error 
between our predicted solution $u_{\text{pred}}^{n}$ and 
the reference solution $u_{\text{ref}}^{n}$, i.e., 
\begin{equation} 
\begin{split}
\| u_{\text{ref}}^{n} - u_{\text{pred}}^{n} \|_{L^2(\Omega)} = &
\left( \int_\Omega \vert u_{\text{ref}}^{n} - u_{\text{pred}}^{n} \vert^2 \, dx \right)^\frac{1}{2}, \\
\| u_{\text{ref}}^{n} - u_{\text{pred}}^{n} \|_{H^1(\Omega)} = &
\left( \int_\Omega \vert \nabla (u_{\text{ref}}^{n} - u_{\text{pred}}^{n}) \vert^2 \, dx \right)^\frac{1}{2}, \\
\end{split} 
\end{equation}

\subsubsection{Experiment 1}

In this experiment, we consider flow in a fixed 
static channelized field $\kappa$ 
and a time-independent source $g$ 
fixed among all the samples. 
The nonlinearity constant is chosen as $\alpha = 20$. 
We use POD to extract 10 dominant modes from $1000$ snapshot solutions 
and construct POD nodal basis functions. 
In the neural network, we simply take the input and output as 
\begin{equation} \mathbf{x} = {c_s}^n \quad \text{ and } \quad \mathbf{y} = {c_s}^{n+1}. \end{equation}
In the generation of samples, 
we consider independent and uniformly distributed initial conditions ${c_s}^0$.
We generate $100$ realizations of initial conditions ${c_s}^0$, 
and evolve the reduced-order dynamic process to obtain ${c_s}^n$ for $n = 1,2,\ldots,10$ . 
We remark that these simulation data provide 
a total of $1000$ samples of the pushforward map. 

We use the $900$ samples given by $90$ realizations as training set and 
the $100$ samples given by $10$ remaining realizations as testing samples. 
Using the training data and a given network architecture, we find a set of optimized parameter $\theta^*$ 
which minimizes the loss function, 
and obtain optimized network parameters $\theta^*$. 
The network $\mathcal{N}$ is then used to predict the $1$-step dynamic, i.e.,
\begin{equation} {c_s}^{n+1} \approx \mathcal{N}({c_s}^n; \theta^*). \end{equation}
We also use the composition of the network $\mathcal{N}$ to predict the final-time solution, i.e.,
\begin{equation} {c_s}^{10} \approx \mathcal{N}(\mathcal{N}(\cdots\mathcal{N}({c_s}^0; \theta^*)\cdots; \theta^*); \theta^*). \end{equation}
We use the same set of training data and testing data 
and compare the performance of different network architectures. 
We examine the performance of the networks by 
the mean of $L^2$ percentage error of the $1$-step prediction and the final-time prediction in the testing samples. 
The error is computed by comparing to the solution formed by the simulation data ${c_s}^n$.

The results are summarized in Table~\ref{tab:exp1.1}. 
It can be observed that if thee network architecture is too simple, 
i.e. contains too few layers or neurons, 
the neural network built may become useless in prediction.

\begin{table}[ht!]
\centering
\begin{tabular}{c|c||c|c}
Layer & Neuron & $1$-step & Final-time \\
\hline
\multirow{3}{*}{3} 
 & 20 & 0.1776 & 3.4501e+09\\
 & 100 & 0.0798 & 6.3276\\
 & 400 & 0.0613 & 4.9855\\
 \hline
\multirow{3}{*}{5} 
 & 20 & 0.1499 & 6.5970e+06\\
 & 100 & 0.0753 & 6.3101\\
 & 400 & 0.0602 & 4.8137\\
 \hline
\multirow{3}{*}{10} 
 & 20 & 0.1024 & 5.4183\\
 & 100 & 0.0750 & 4.3271\\
 & 400 & 0.0609 & 1.8834\\
 \end{tabular}
\caption{Mean of $L^2$ percentage error with different network architectures in Experiment 1}
\label{tab:exp1.1}
\end{table}

\subsubsection{Experiment 2}

\rev{
In the second experiment, we consider flow in 
static channelized fields $\kappa$ 
and a time-independent source $g$
fixed among all the samples. 
The nonlinearity constant is chosen as $\alpha = 10$. 
The coefficient fields $\kappa$ 
differ in the conductivity value in channels. 
The high conductivity values in the two channels are parametrized by 
\begin{equation}
\begin{split}
\kappa_1 & = 10000 e^{\eta_1}, \\
\kappa_2 & = 5000 e^{\eta_2}, 
\end{split}
\end{equation}
where 
$\eta = (\eta_1, \eta_2)$ is taken from a uniform distribution in $[-0.5, 0.5]^2$. 
We use POD to extract 5 dominant modes from $1000$ snapshot solutions 
and construct POD nodal basis functions. 
In the neural network, we simply take the input and output as 
\begin{equation}
\mathbf{x} = ({c_s}^n, \eta) \quad \text{ and } \quad \mathbf{y} = {c_s}^{n+1}.
\end{equation}
We generate $100$ realizations of initial conditions ${c_s}^0$ and parameters $\eta$, 
and evolve the reduced-order dynamic process to obtain ${c_s}^n$ for $n = 1,2,\ldots,10$. 
We remark that these simulation data provide 
a total of $1000$ samples of the pushforward map. 
We use the $900$ samples given by $90$ realizations as training set and 
the $100$ samples given by $10$ remaining realizations as testing samples. 
Using the training data and a given network architecture, 
we find a set of optimized parameter $\theta^*$ 
which minimizes the loss function, 
and obtain optimized network parameters $\theta^*$. 
The network $\mathcal{N}$ is then used to predict the $1$-step dynamic, i.e.,
\begin{equation}
{c_s}^{n+1} \approx \mathcal{N}({c_s}^n, \eta; \theta^*).
\end{equation}
We also use the composition of the network $\mathcal{N}$ to predict the final-time solution, i.e.,
\begin{equation}
{c_s}^{10} \approx \mathcal{N}(\mathcal{N}(\cdots\mathcal{N}({c_s}^0,\eta; \theta^*)\cdots,\eta; \theta^*),\eta; \theta^*).
\end{equation}
}

\rev{
In this example, we investigate the advantage of our approach of
combining deep learning with POD nodal basis functions. 
Instead of using the coefficients of the solution 
with respect to POD nodal basis functions $\{ \psi_k \}_{k=1}^m$ for representing the flow dynamics, 
one can also use other discretizations,
for example, the standard bilinear elements nodal functions 
or the POD basis functions $\{ \mathbf{v}_j \}_{j=1}^m$ 
Using the same idea as in Section~\ref{sec:method}, 
we can learn from the respective data and 
construct corresponding neural networks for approximations. 
In this experiment, we compare the training cost and the performance of the neural networks 
using different underlying discretizations, 
by using the same set of training data and testing data. 
All networks consist of 3 hidden layers of 20 neurons 
and are trained in 500 epochs. 
We examine the performance of the networks by 
comparing the $1$-step prediction and 
the final-time prediction to the corresponding numerical method. 
}

\rev{
A comparison of discretizations is presented in Table~\ref{tab:exp1.2}, which 
suggest that the model reduction technique brings several advantages 
to neural network approximation of flow dynamics. 
First, the use of POD reduces the number of 
trainable parameters in the network 
and thus shortening the elapsed time for network training. 
In our simple experiment, as shown in Table~\ref{tab:exp1.2}, 
the elapsed time for training the networks in POD reduced-order models
is around 1/10 of elapsed time for training the networks in the standard nodal coordinates. 
Second, instead of extracting features solely in the learning process, 
the reduced order model predefines some features which 
are important in representing the flow and facilitates 
the learning process. This allows the information propagates more 
easily through the multi-layer networks and provides a smaller prediction error. 
Lastly, learning the evolution in the standard nodal coordinates becomes 
infeasible in large-scale computation. 
Both elapsed runtime for sample generation 
and memory required for sample storage 
grow dramatically with increased number of degree of freedom. 
The reduced-order model provides a cheap alternative 
for learning the flow dynamics in this scenario. 
\rev{As shown in Table~\ref{tab:exp1.2}, the CPU time for one forward run in the full model 
is 0.4499 seconds, which is short due to the simplicity of the linearization scheme in the simple experiment. 
However, with the reduced order model, the CPU time for a single forward run 
is reduced to 0.0003 seconds. 
We remark that the use of reduced-order models will be even more advantageous 
in complicated problems. 
For example, for repeatedly modelling highly nonlinear flows in highly heterogeneous flows, 
the nonlinear solver in the high-fidelity space will be computationally expensive. }
Moreover, the prediction error using POD nodal basis functions $\{ \psi_k \}_{k=1}^m$
is smaller than using the original POD basis functions $\{ \mathbf{v}_j \}_{j=1}^m$.
This suggests that nodal values provide a more stable and
well-conditioned coordinate system. 
}

\rev{
\begin{table}[ht!]
\centering
\begin{tabular}{c||c|c|c}
Coordinates & Standard nodal & POD & POD nodal \\ 
\hline \hline
Dimension & 9801 & 5 & 5 \\
\rev{Forward runtime (seconds)} & \rev{0.4499} & \rev{0.0003} & \rev{0.0003} \\
\hline
\# trainable parameters & 403161 & 1525 & 1525 \\
\rev{Training} time (seconds) & 587.14 & 62.60 & 57.56 \\
\hline
$L^2$ error for $1$-step & 0.9529\% & 0.5751\% & 0.3957\% \\
$H^1$ error for $1$-step & 2.9588\% & 0.6020\% & 0.4395\% \\
\hline
$L^2$ error for final time & 4.8563\% & 3.4943\% & 3.0266\% \\
$H^1$ error for final time & 5.9270\% & 3.7307\% & 3.3762\% \\
\end{tabular}
\caption{History of training cost and prediction error with different discretization in Experiment 2.} 
\label{tab:exp1.2}
\end{table}
}

\rev{
\begin{remark}
The wide neural network using standard nodal coordinates 
can be viewed as a generalization of dynamic mode decomposition (DMD) \cite{schmid2010dmd}. 
DMD is a dimensionality reduction technique which extracting dynamical features from flow data. 
Given a sequence of snapshots $\{ u_h^0, u_h^1, \ldots, u_h^K \}$, 
DMD seeks a linear mapping $A$ which fits the snapshots by $u_h^{n+1} = A u_h^n$, 
which can be seen as the simplest neural network with 
linear activation function and without bias and hidden layers that maps $u_h^n$ to $u_h^{n+1}$. 
Optimal mode decomposition (OMD) \cite{wynn2014omd}, a variant of DMD, 
seeks a linear mapping $A$ with a user-defined rank $k$, which is equivalent to 
seek a wide neural network in the form 
\begin{equation}
u_h^{n+1} = W_2 W_1 u_h^n,
\end{equation}
i.e. a 2-layer network with linear activation and no bias, and with $k$ neurons in the 
immediate hidden layer. In this sense, we can build more general neural networks 
than DMD or OMD, which provides higher interpretability for more complex 
and nonlinear flow dynamics. 
\end{remark}
}

\subsection{Application to two-phase flow}

\rev{
In this example, we consider the evolution of water saturation a two-phase oil-water flow in porous media. 
Following \cite{poddeim15,efendiev2016online}, we employ POD for model order reduction in the unknowns and 
Discrete Empirical Interpolation Method (DEIM) for approximations of nonlinear source functions. 
Let $\Omega$ be the reservoir domain, in which 
the two phases, namely water (denoted by subscript $w$) and 
oil (denoted by subscript $o$), of the flow are immiscible. 
The flow equation of the phase velocity $u_a$ is described by the Darcy's law:
\begin{equation}
u_a = -\mu_a^{-1} k_{ra}(s) \kappa \nabla p, 
\label{eq:darcy}
\end{equation}
where $\kappa$ is the permeability tensor, 
$k_{ra}$ is the relative permeability to phase $a$, 
$s$ is the (water) saturation, 
and $p$ is the pressure, 
for each phase $a = w,o$. 
Here, we assume the displacement is dominated by viscous effects, 
so that the gravitational acceleration, capillary pressure effects 
and compressibility effects can be neglected. 
The Darcy's law is coupled with the conservation of mass and 
this yields a initial-boundary value problem, 
known as the pressure-saturation equations: 
\begin{equation}
\begin{split}
-\nabla \cdot (\lambda(s) \kappa \nabla p) & = q_w + q_o \text{ in } \Omega, \\
\phi \dfrac{\partial s}{\partial t} + \nabla \cdot (f_w(s) u) & = \dfrac{q_w}{p_w} \text{ in } \Omega, \\
u \cdot n & = 0 \text{ on } \partial \Omega, \\
s & = 0 \text{ at } t = 0.
\end{split}
\label{eq:two-phase}
\end{equation}
Here $q_w$ and $q_o$ are the volumetric source function 
acting to the water and oil respectively, 
$\phi$ is the porosity, 
$u = u_w + u_o$ is the total velocity, 
$\lambda$ is the total mobility defined by 
\begin{equation}
\lambda(s) = \lambda_w(s) + \lambda_o(s) = \dfrac{k_{rw}(s)}{\mu_w} + \dfrac{k_{ro}(s)}{\mu_o},
\end{equation}
and $f_w$ is the flux defined by 
\begin{equation}
f_w(s) = \dfrac{\lambda_w(s)}{\lambda(s)} = \dfrac{k_{rw}(s)}{k_{rw}(s) + \frac{\mu_w}{\mu_o} k_{ro}(s)}. 
\end{equation}
To solve \eqref{eq:two-phase} numerically, in each time step, 
we sequentially update pressure, saturation and the physical quantities. 
First, we solve for the pressure $p^{n+1}$ and the velocity $u^{n+1}$ 
with a mixed finite element method. 
For details of the mixed finite element method, 
we refer the readers to \cite{mixedfem,poddeim15}. 
Next, we solve for the saturation with a mass conservative finite volume discretization. 
In a cell $\Omega_i$, we write the cell average of the water saturation as $s_i^{n+1}$ 
and the common edges of the cells $\Omega_i$ and $\Omega_j$ as $\gamma_{ij}$.  
Then we approximate the flux over $\gamma_{ij}$ by the upstream weighting
\begin{equation}
F_{ij}(s^{n+1}) \approx \int_{\gamma_{ij}} 
f_w(s_i^{n+1}) (u^{n+1} \cdot n_{ij})^+ + 
f_w(s_j^{n+1}) (u^{n+1} \cdot n_{ij})^-, 
\end{equation}
where $n_{ij}$ is the outward normal vector on $\gamma_{ij}$ and 
the superscripts $+$ and $-$ denote the positive part and the negative part respectively. 
The evolution of saturation is given by the nonlinear equation
\begin{equation}
s_i^{n+1} = s_i^n + \dfrac{\Delta t}{\vert \Omega_i \vert} 
\left( q_i^+ - \sum_j F_{ij}(s^{n+1}) + f_w(s_i^{n+1}) q_i^- \right), 
\end{equation}
which is numerically solved by Newton iterations. 
The map from $s^n$ to $s^{n+1}$ is therefore very complex and nonlinear. 
Using POD global model reduction, 
we obtain POD basis $(\mathbf{V}_p, \mathbf{V}_u, \mathbf{V}_s)$ 
for pressure, velocity and saturation respectively. 
Here we briefly mention that the nonlinear functions $s \mapsto (\lambda(s) \kappa)^{-1}$ 
and $s \mapsto f_w(s)$ are evaluated in a projected subspace 
using Discrete Empirical Interpolation Method (DEIM). 
For details of DEIM, the reader is referred to \cite{poddeim15,efendiev2016online}. 
Using the reduced-order model, forward simulations are performed 
and used as samples in the neural network. 
}

\rev{
In our simulations, the relative permeability to both phases are quadratic in $s$.
The flow is simulated for 300 days and the saturation information is collected and saved every 30 days.
Realizations of flow dynamics are computed using 
independent and uniformly distributed initial conditions, permeability field in well rates. 
\rev{The dimension of the finite element space is 2500.}
We use POD to extract 25 dominant modes for saturation 
and construct POD nodal basis functions. 
Simulation data of the dynamic process under the reduced-order model 
are then obtained and used in the training set. 
In some cases, observation data are supplemented in the training set. 
Different forms of inputs, depending on situations, are investigated. 
Using these data as samples, universal multi-layer networks are trained to 
approximate the flow dynamics. 
We use the trained networks to predict the output with some 
new unseen inputs, and reconstruct the numerical solution 
using the predicted coefficients. 
We examine the quality of our networks by computing the $L^2$ error 
between our predicted solution $s_{\text{pred}}^{n}$ and 
the reference solution $s_{\text{ref}}^{n}$, i.e., 
\begin{equation} 
\begin{split}
\| s_{\text{ref}}^{n} - s_{\text{pred}}^{n} \|_{L^2(\Omega)} = &
\left( \int_\Omega \vert s_{\text{ref}}^{n} - s_{\text{pred}}^{n} \vert^2 \, dx \right)^\frac{1}{2}.
\end{split} 
\end{equation}
}

\subsubsection{Experiment 3}

\begin{figure}[!h]
	\centering
	\begin{subfigure}{.32\textwidth}
		\centering
		\includegraphics[width=.9\linewidth]{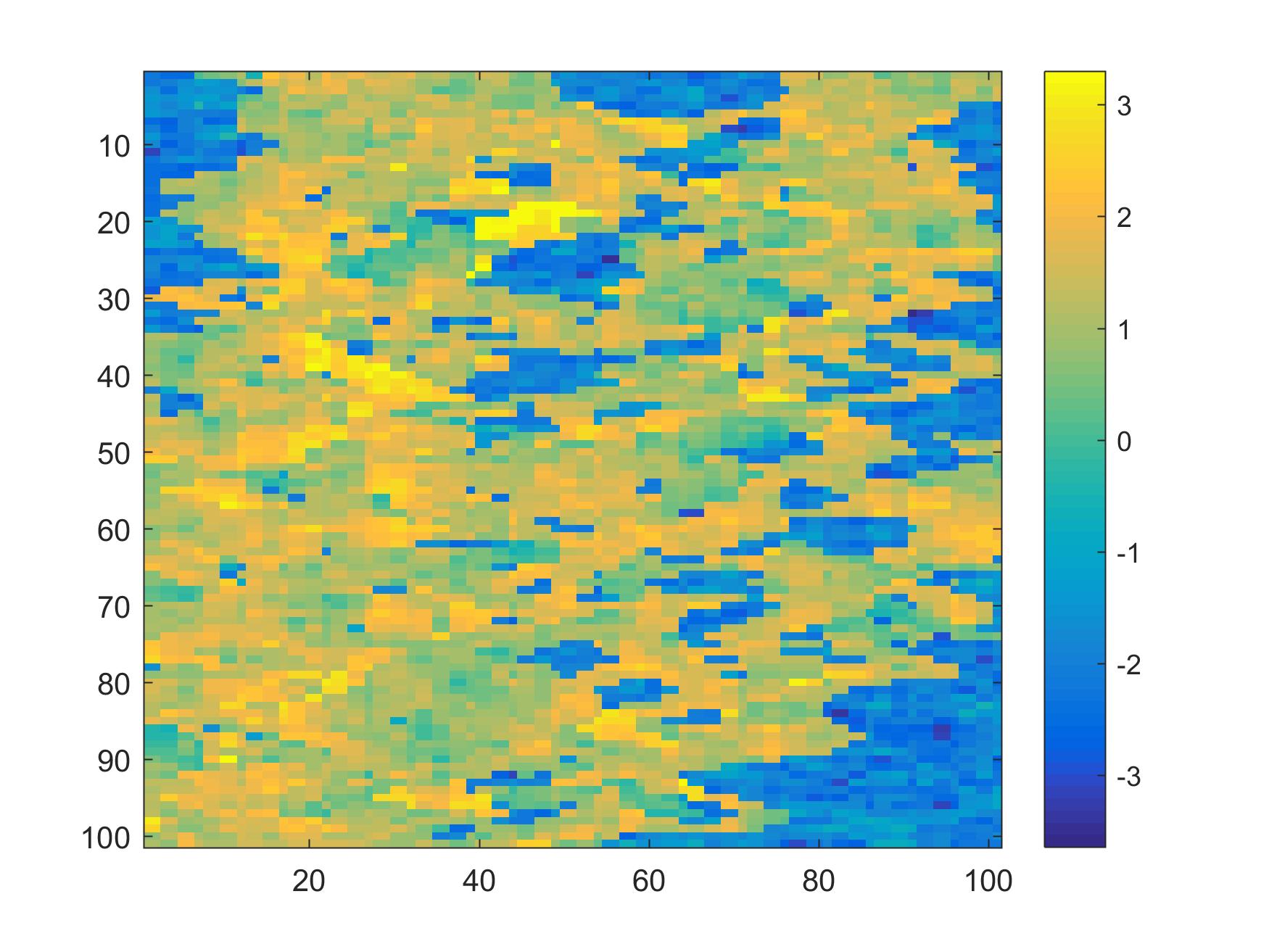}
		\caption{$\log(\kappa)$ for $\eta = 0$}	
	\end{subfigure}
	\begin{subfigure}{.32\textwidth}
		\centering
		\includegraphics[width=.9\linewidth]{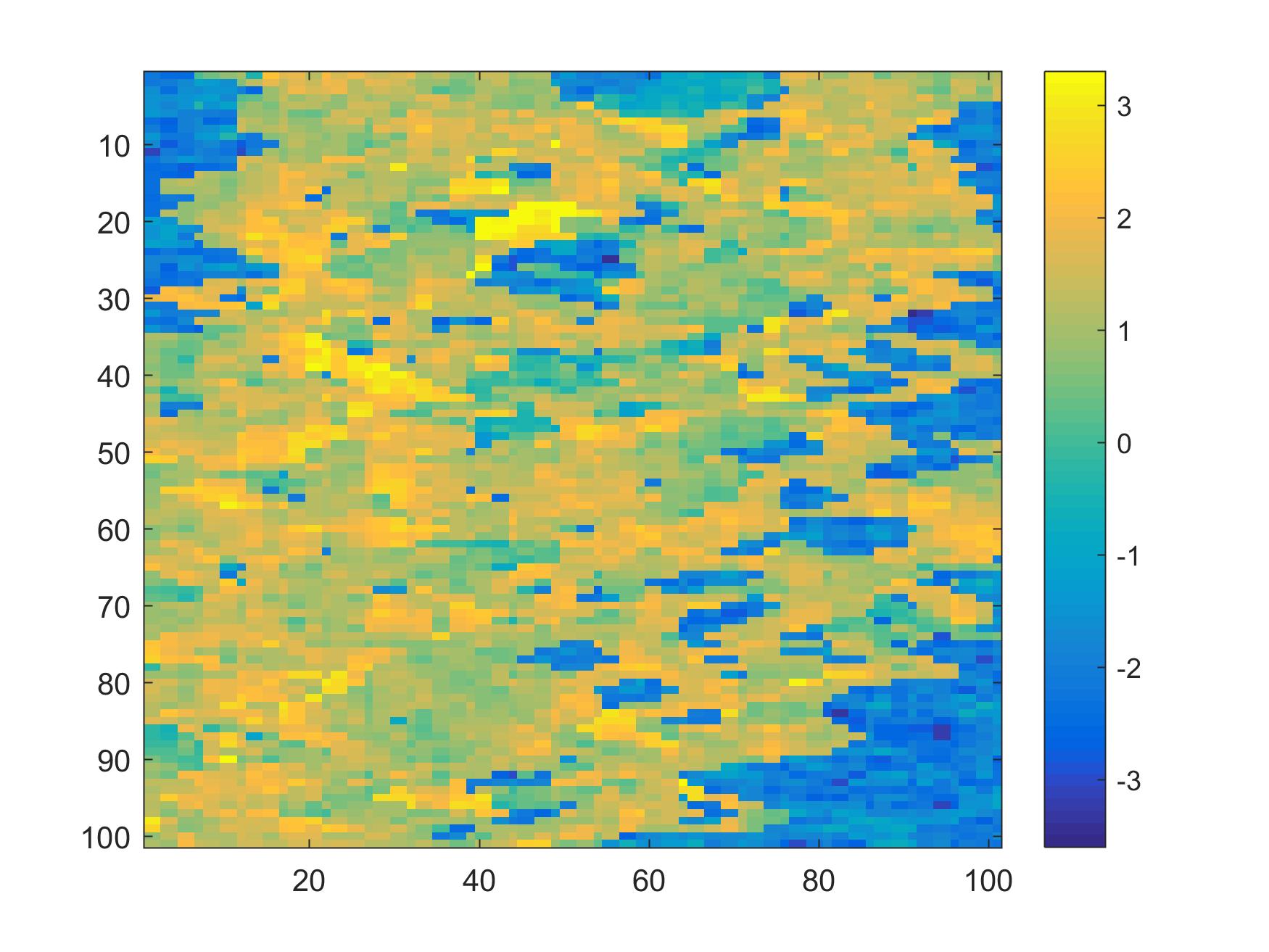}
		\caption{$\log(\kappa)$ for $\eta = 0.05$}
	\end{subfigure}
	\begin{subfigure}{.32\textwidth}
		\centering
		\includegraphics[width=.9\linewidth]{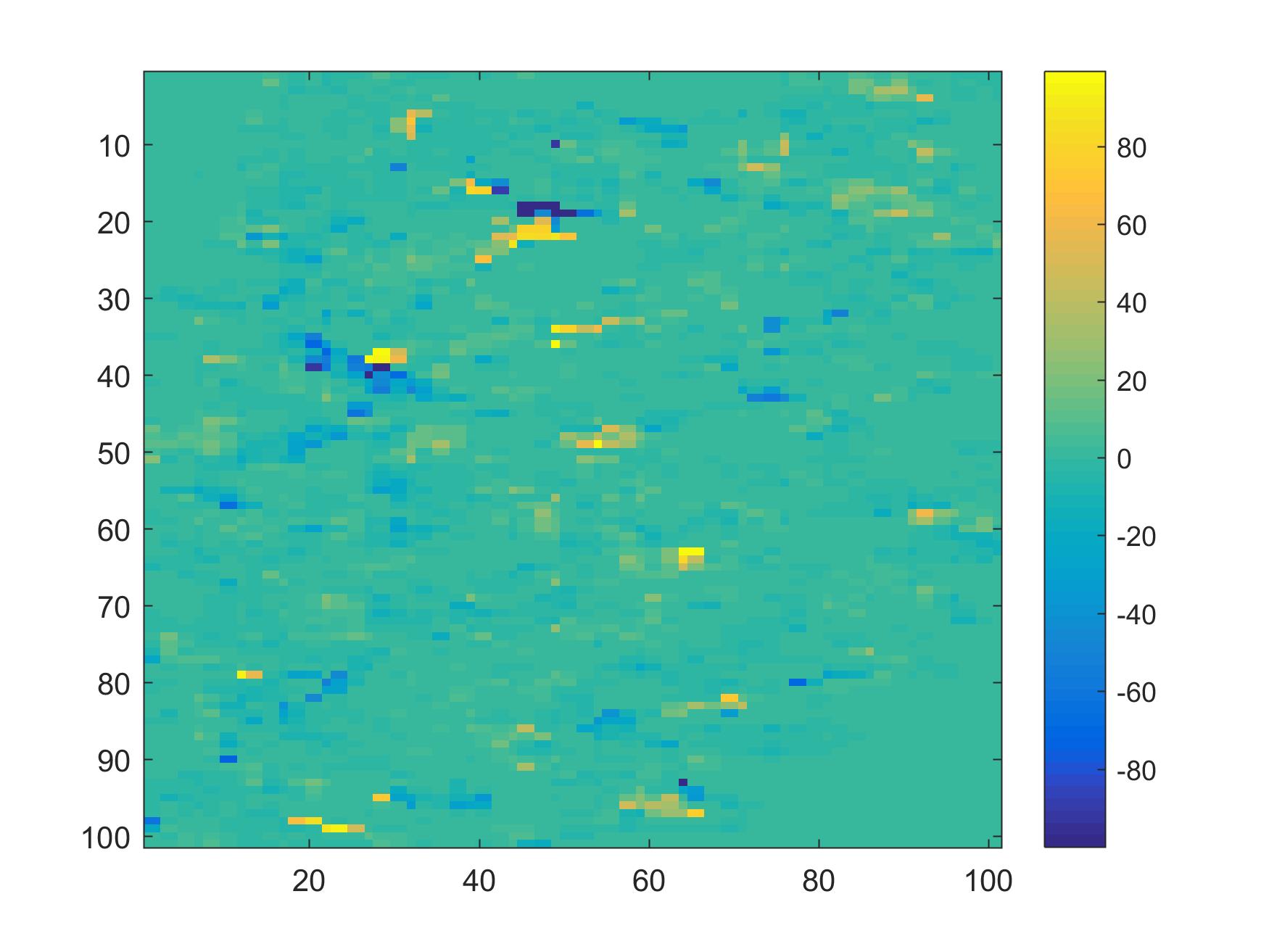}
		\caption{Difference of $\kappa$}
	\end{subfigure}
	\caption{Figures of permeability fields used in two-phase flow.}
	\label{fig:spe10}
\end{figure}

\rev{
In this experiment, we consider flow in varying permeability field $\kappa$ 
and with varying source function. 
We consider the case with uncertainties in the permeability fields, 
which can be parametrized by a convex combination
\begin{equation}
\kappa = \eta \kappa_1 + (1-\eta) \kappa_2, 
\end{equation}
where the parameter $\eta$ is taken from a uniform distribution in $[0,0.1]$, 
and the fixed permeability fields $\kappa_1$ and $\kappa_2$ 
are taken from SPE10 benchmarks. 
Samples for permeability fields used in Experiment 3 
are depicted in Figure~\ref{fig:spe10}. 
The source function $q^{n+1}$ is external forcing at certain well positions. 
There are 1 injection well and altogether 4 production wells. 
In each simulation, only 1 production well is active. 
The source function $q^{n+1}$ can therefore be parametrized by the well rates in $\mathbb{R}^4$. 
In the neural network, we simply take the input and output as 
\begin{equation} \mathbf{x} = ({c_s}^n,\eta,q^{n+1}) \quad \text{ and } \quad \mathbf{y} = {c_s}^{n+1}. \end{equation}
}

\rev{
In this example, we investigate the data-driven neural networks. 
We assume different availability of observation data, 
which is synthetic data from expensive fine-scale calculations. 
\rev{The CPU time for a single forward run of fine-scale observation data is 49.5929 seconds, 
while that of POD simulation data is 17.9847 seconds.}
In practice, observation data can be obtained from physical measurement. 
We compare the performance of observation network $\mathcal{N}_o$, 
mixed network $\mathcal{N}_m$ and simulation network $\mathcal{N}_s$ 
by taking the fine-scale solution as reference. 
The observational data therefore serves as a corrector of the neural network
in the scenario of using simulation data supplemented by limited observation data. 
All the networks are trained with $3520$ samples and 
identical neural network architectures, where each intermediate hidden layer 
consists of 300 neurons. 
For the mixed network, one half of training data is from simulation 
and the other half is from observation. 
Using the training data, we find a set of optimized parameter $\theta^*$ 
which minimizes the loss function, 
and obtain optimized network parameters $\theta^*$. 
The network $\mathcal{N}$ is then used to predict the $1$-step dynamic, i.e.,
\begin{equation} {c_s}^{n+1} \approx \mathcal{N}({c_s}^n, \eta, q^{n+1}; \theta^*). \end{equation}
We also use the composition of the network $\mathcal{N}$ to predict the final-time solution, i.e.,
\begin{equation} {c_s}^{10} \approx \mathcal{N}(\mathcal{N}(\cdots\mathcal{N}({c_s}^0,\eta,q^1; \theta^*)\cdots,\eta,q^{9}; \theta^*),\eta,q^{10}; \theta^*). \end{equation}
}

\rev{
The results are summarized in Tables~\ref{tab:exp3}. It can be 
observed that the prediction of 1-step dynamics is not much affected 
by the network architecture or the type of training data. 
However, the observation data provides improvement to the quality 
of prediction of final-time dynamics. 
On the other hand, the return of increasing number of layers is eventually marginal,
but shallow networks are not able to fully interpret all the intrinsic data. 
}

\begin{table}[ht!]
\centering
\begin{tabular}{c||c|c|c||c|c|c}
& \multicolumn{3}{c||}{1-step} & \multicolumn{3}{c}{Final time} \\
\cline{2-4} \cline{5-7}
Layer & Simulation & Mixed & Observation & Simulation & Mixed & Observation \\
\hline \hline
1 & 1.15\% & 2.25\% & 1.77\% & 18.98\% & 11.49\% & 2.77\% \\
3 & 0.59\% & 1.51\% & 1.33\% & 14.97\% & 8.14\% & 1.65\% \\
5 & 0.78\% & 1.33\% & 1.05\% & 15.00\% & 8.07\% & 1.06\% \\ 
10 & 0.80\% & 1.26\% & 1.52\% & 14.80\% & 8.13\% & 1.66\%
\end{tabular}
\caption{Mean of $L^2$ percentage error with different training data in Experiment 3.}
\label{tab:exp3}
\end{table}

\section{Conclusion}

In this paper, we combine some POD techniques with deep learning concepts
in the simulations for flows in porous media. The observation data 
is given at some locations. We construct POD modes such that 
the degrees of freedom represent the values of the solution at centain 
locations. Furthermore, we write the solution at the current time
as a multi-layer network that depends on the solution at the initial
time and input parameters, such as well rates and permeability fields.
This provides a natural framework for applying deep learning techniques
 for flows in channelized media. We provide the details of our method
and present numerical results. In all numerical results, we study
nonlinear flow equation in channelized media and consider various channel
configurations. Our results show that multi-layer network provides
an accurate approximation of the forward map and can incorporate the observed
data. Moreover, by incorporating some observed data (from true model)
and some computational data, we modify the reduced-order model.
This way, one can use the observed data to modify reduced-order
models which honor the observed data.

\bibliographystyle{plain}
\bibliography{references}

\begin{thebibliography}{10}

\bibitem{nonlinear_AM2015}
Manal Alotaibi, Victor~M. Calo, Yalchin Efendiev, Juan Galvis, and Mehdi
  Ghommem.
\newblock Global-local nonlinear model reduction for flows in heterogeneous
  porous media.
\newblock {\em Computer Methods in Applied Mechanics and Engineering},
  292:122--137, 2015.

\bibitem{mixedfem}
D.~Boffi, F.~Brezzi, and M.~Fortin.
\newblock {\em Mixed Finite Element Methods and Applications}.
\newblock Springer-Verlag, Heidelberg, 2013.

\bibitem{cd10}
M.~Cardoso and L.~Durlofsky.
\newblock Linearized reduced-order models for subsurface flow simulation.
\newblock {\em Journal of Computational Physics}, 229, 2010.

\bibitem{cardoso2009development}
MA~Cardoso, LJ~Durlofsky, and P~Sarma.
\newblock Development and application of reduced-order modeling procedures for
  subsurface flow simulation.
\newblock {\em International journal for numerical methods in engineering},
  77(9):1322--1350, 2009.

\bibitem{celia1990general}
Michael~A Celia, Efthimios~T Bouloutas, and Rebecca~L Zarba.
\newblock A general mass-conservative numerical solution for the unsaturated
  flow equation.
\newblock {\em Water resources research}, 26(7):1483--1496, 1990.

\bibitem{chollet2015keras}
Fran\c{c}ois Chollet et~al.
\newblock Keras.
\newblock \url{https://keras.io}, 2015.

\bibitem{chung2015residual}
E.~T. Chung, Y.~Efendiev, and W.~T. Leung.
\newblock Residual-driven online generalized multiscale finite element methods.
\newblock {\em Journal of Computational Physics}, 302:176--190, 2015.

\bibitem{Csaji2001}
Balázs~Csanád Csáji.
\newblock Approximation with artificial neural networks.
\newblock {\em Faculty of Sciences, Etvs Lornd University}, 24(48), 2001.

\bibitem{Cybenko1989}
G.~Cybenko.
\newblock Approximations by superpositions of sigmoidal functions.
\newblock {\em Mathematics of Control, Signals, and Systems}, 2(4):303--314,
  1989.

\bibitem{dostert2009efficient}
P~Dostert, Y~Efendiev, and B~Mohanty.
\newblock Efficient uncertainty quantification techniques in inverse problems
  for richards’ equation using coarse-scale simulation models.
\newblock {\em Advances in water resources}, 32(3):329--339, 2009.

\bibitem{E_deepRitz}
Weinan E and Bing Yu.
\newblock The deep {R}itz method: A deep learning-based numerical algorithm for
  solving variational problems.
\newblock {\em Communications in Mathematics and Statistics}, 6(1):1--12, 2018.

\bibitem{efendiev2005efficient}
Y~Efendiev, A~Datta-Gupta, V~Ginting, X~Ma, and B~Mallick.
\newblock An efficient two-stage markov chain monte carlo method for dynamic
  data integration.
\newblock {\em Water Resources Research}, 41(12), 2005.

\bibitem{efendiev2012local}
Yalchin Efendiev, Juan Galvis, and Eduardo Gildin.
\newblock Local--global multiscale model reduction for flows in high-contrast
  heterogeneous media.
\newblock {\em Journal of Computational Physics}, 231(24):8100--8113, 2012.

\bibitem{efendiev2016online}
Yalchin Efendiev, Eduardo Gildin, and Yanfang Yang.
\newblock Online adaptive local-global model reduction for flows in
  heterogeneous porous media.
\newblock {\em Computation}, 4(2):22, 2016.

\bibitem{gardner1958some}
WR~Gardner.
\newblock Some steady-state solutions of the unsaturated moisture flow equation
  with application to evaporation from a water table.
\newblock {\em Soil science}, 85(4):228--232, 1958.

\bibitem{global14}
Mehdi Ghommem, Victor~M. Calo, and Yalchin Efendiev.
\newblock Mode decomposition methods for flows in high-contrast porous media. a
  global approach.
\newblock {\em Journal of Computational Physics}, 257:400--413, 2014.

\bibitem{globallocal13}
Mehdi Ghommem, Michael Presho, Victor~M. Calo, and Yalchin Efendiev.
\newblock Mode decomposition methods for flows in high-contrast porous media.
  global-local approach.
\newblock {\em Journal of Computational Physics}, 253:226--238, 2013.

\bibitem{glorot11}
Xavier Glorot, Antoine Bordes, and Yoshua Bengio.
\newblock Deep sparse rectifier neural networks.
\newblock In {\em Proceedings of the Fourteenth International Conference on
  Artificial Intelligence and Statistics}, pages 315--323. PMLR, 2011.

\bibitem{goodfellow2016deep}
Ian Goodfellow, Yoshua Bengio, Aaron Courville, and Yoshua Bengio.
\newblock {\em Deep learning}, volume~1.
\newblock MIT press Cambridge, 2016.

\bibitem{Hanin2017}
Boris Hanin.
\newblock Universal function approximation by deep neural nets with bounded
  width and relu activations.
\newblock {\em arXiv:1708.02691}, 2017.

\bibitem{volkwein05}
M.~Hinze and S.~Volkwein.
\newblock Proper orthogonal decomposition surrogate models for nonlinear
  dynamical systems: error estimates and suboptimal control.
\newblock In P.~Benner, V.~Mehrmann, and D.C. Sorensen, editors, {\em Dimension
  Reduction of Large-Scale Systems}, volume~45 of {\em Lecture Notes in
  Computational Science and Engineering}, pages 261--306. Springer Berlin
  Heidelberg, 2005.

\bibitem{Hornik1991}
Kurt Hornik.
\newblock Approximation capabilities of multilayer feedforward networks.
\newblock {\em Neural Networks}, 4(2):251–257, 1991.

\bibitem{jansen2017use}
Jan~Dirk Jansen and Louis~J Durlofsky.
\newblock Use of reduced-order models in well control optimization.
\newblock {\em Optimization and Engineering}, 18(1):105--132, 2017.

\bibitem{Kerschen2005}
Gaetan Kerschen, Jean-claude Golinval, Alexander~F. Vakakis, and Lawrence~A.
  Bergman.
\newblock The method of proper orthogonal decomposition for dynamical
  characterization and order reduction of mechanical systems: An overview.
\newblock {\em Nonlinear Dynamics}, 41(1):147--169, 2005.

\bibitem{Ying_paraPDE}
Yuehaw Khoo, Jianfeng Lu, and Lexing Ying.
\newblock Solving parametric pde problems with artificial neural networks.
\newblock {\em arXiv:1707.03351}, 2017.

\bibitem{adam}
Diederik~P. Kingma and Jimmy Ba.
\newblock Adam: A method for stochastic optimization.
\newblock {\em arXiv preprint arXiv:1412.6980}, 2014.

\bibitem{lecun2015deep}
Yann LeCun, Yoshua Bengio, and Geoffrey Hinton.
\newblock Deep learning.
\newblock {\em nature}, 521(7553):436, 2015.

\bibitem{Shi_resnet}
Zhen Li and Zuoqiang Shi.
\newblock Deep residual learning and pdes on manifold.
\newblock {\em arXiv:1708.05115.}, 2017.

\bibitem{relu}
A.L. Maas, A.Y. Hannun, and A.Y. Ng.
\newblock Rectifier nonlinearities improve neural network acoustic models.
\newblock {\em Proc. icml}, 30(1):3, 2013.

\bibitem{Poggio2016}
Hrushikesh Mhaskar, Qianli Liao, and Tomaso Poggio.
\newblock Learning functions: when is deep better than shallow.
\newblock {\em arXiv preprint arXiv:1603.00988}, 2016.

\bibitem{richards1931capillary}
Lorenzo~Adolph Richards.
\newblock Capillary conduction of liquids through porous mediums.
\newblock {\em physics}, 1(5):318--333, 1931.

\bibitem{schmid2010dmd}
P.~J. Schmid.
\newblock Dynamic mode decomposition of numerical and experimental data.
\newblock {\em Journal of Fluid Mechanics}, 656:5--28, 2010.

\bibitem{schmidhuber2015deep}
J{\"u}rgen Schmidhuber.
\newblock Deep learning in neural networks: An overview.
\newblock {\em Neural networks}, 61:85--117, 2015.

\bibitem{Telgrasky2016}
M.~Telgrasky.
\newblock Benefits of depth in neural nets.
\newblock {\em JMLR: Workshop and Conference Proceedings}, 49(123), 2016.

\bibitem{trehan2016trajectory}
Sumeet Trehan and Louis~J Durlofsky.
\newblock Trajectory piecewise quadratic reduced-order model for subsurface
  flow, with application to pde-constrained optimization.
\newblock {\em Journal of Computational Physics}, 326:446--473, 2016.

\bibitem{van2006reduced}
Jorn~FM van Doren, Renato Markovinovi{\'c}, and Jan-Dirk Jansen.
\newblock Reduced-order optimal control of water flooding using proper
  orthogonal decomposition.
\newblock {\em Computational Geosciences}, 10(1):137--158, 2006.

\bibitem{van1980closed}
M~Th Van~Genuchten.
\newblock A closed-form equation for predicting the hydraulic conductivity of
  unsaturated soils 1.
\newblock {\em Soil science society of America journal}, 44(5):892--898, 1980.

\bibitem{vo2014new}
Hai~X Vo and Louis~J Durlofsky.
\newblock A new differentiable parameterization based on principal component
  analysis for the low-dimensional representation of complex geological models.
\newblock {\em Mathematical Geosciences}, 46(7):775--813, 2014.

\bibitem{wynn2014omd}
A.~Wynn, D.~S. Pearson, B.~Ganapathisubramani, and P.~J. Goulart.
\newblock Optimal mode decomposition for unsteady flows.
\newblock {\em Journal of Fluid Mechanics}, 733:473--503, 2013.

\bibitem{poddeim15}
Yanfang Yang, Mohammadreza Ghasemi, Eduardo Gildin, Yalchin Efendiev, and
  Victor Calo.
\newblock Fast multiscale reservoir simulations with pod-deim model reduction.
\newblock {\em SPE Journal}, 21(06):2141--2154, 2016.

\end{thebibliography}

\end{document}